\documentclass[a4paper]{article}
\usepackage{amsmath,amssymb,amsthm,amsfonts,amscd}
\usepackage{geometry}

\title{The cohomology of lattices in ${\rm SL}(2,\mathbb{C})$}
\author{
Tobias Finis \thanks{e-mail: {\tt finis@math.uni-duesseldorf.de}}\\
Mathematisches Institut\\
Heinrich-Heine-Universit\"at\\
D-40225 D\"usseldorf
\and 
Fritz Grunewald \thanks{e-mail: {\tt fritz@math.uni-duesseldorf.de}}\\
Mathematisches Institut\\
Heinrich-Heine-Universit\"at\\
D-40225 D\"usseldorf
\and
Paulo Tirao \thanks{e-mail: {\tt ptirao@famaf.unc.edu.ar}}  \\
Facultad de Matem\'atica, Astronom\'ia y F\'isica\\
Universidad Nacional de C\'ordoba\\
Ciudad Universitaria, C\'ordoba, Argentina
}     

\date{August 8, 2008}

\numberwithin{equation}{section} 
 
\theoremstyle{plain} 
\newtheorem{theorem} {Theorem} [section]
\newtheorem{lemma} [theorem] {Lemma}
\newtheorem{corollary} [theorem] {Corollary}
\newtheorem{proposition} [theorem] {Proposition}

\newtheorem{question} [theorem] {Question}

\newtheorem{definition} [theorem]{Definition}

\newtheorem{example} [theorem] {Example}

\renewcommand{\b}{\mathfrak{b}}
\newcommand{\prim}{\mathfrak{p}}
\newcommand \A {\mathbb A}
\newcommand \C {\mathbb C}

\newcommand \Q {\mathbb Q}
\newcommand \Z {\mathbb Z}
\newcommand \F {\mathbb F}
\newcommand \N {\mathbb N}
\newcommand \PP {\mathbb P}
\newcommand{\Ram}{{\cal R}}
\newcommand \G {\Gamma}

\renewcommand \a {\mathfrak{a}}

\renewcommand \O {\mathcal{O}}
\renewcommand \o {\mathcal{O}}

\newcommand \Sym {\operatorname{Sym}}
\newcommand \Der {\operatorname{Der}}
\newcommand \IDer {\operatorname{IDer}}
\newcommand \PSL {\operatorname{PSL}}
\newcommand \SL {\operatorname{SL}}
\newcommand \GL {\operatorname{GL}}

\newcommand \tr {\operatorname{tr}}
\newcommand \res {\operatorname{res}}

\newcommand \w {\omega}

\newcommand{\Hom}{\operatorname{Hom}}
\newcommand{\Autf}{{\cal A}}

\begin{document}\maketitle

\begin{abstract} 
This paper contains both theoretical results and experimental data 
on the behavior of the dimensions of the cohomology spaces
$H^1(\G,E_n)$, where $\Gamma$ is a lattice in $\SL(2,\C)$ and 
$E_n = \Sym^n\otimes \overline{\Sym}{}^n$, 
$n\in \N\cup \{0\}$, is one of the standard self-dual modules. In the case
$\Gamma = \SL(2,\O)$ for the ring of integers $\O$ in an imaginary
quadratic number field, we make the theory of lifting explicit and obtain 
lower bounds linear in $n$. We have accumulated a large amount of experimental data 
in this case, as well as for some geometrically constructed and mostly non-arithmetic groups. The computations for $\SL(2,\O)$ lead us to discover two instances with non-lifted classes in the cohomology.  
We also derive an upper bound of size $O(n^2 / \log n)$ for any fixed lattice $\G$ in the general case. We discuss a number of new questions and conjectures suggested by our results and our experimental data.
\newline 

\end{abstract}

\tableofcontents

\section{Introduction}\label{intro}

For a semisimple Lie group $G$ and a lattice $\Gamma$ in $G$ 
(i.e. a discrete subgroup of finite covolume), it is natural to 
consider the cohomology groups $H^* (\Gamma, E)$ of $\Gamma$ with 
coefficients in finite-dimensional representation spaces $E$ of $G$. 
If $K$ is a maximal compact subgroup of $G$ and  $X = G / K$ the 
associated Riemannian symmetric space, these cohomology groups are 
canonically isomorphic to the cohomology groups of the quotient 
orbifold $\Gamma \backslash X$ with coefficients in the 
local system associated to $E$. By the main result of 
\cite{FrankeSchwermer}, at least for arithmetic $\Gamma$
the cohomology of $\Gamma$ can be described by automorphic forms. 
For the contribution of the cuspidal spectrum one has (for general $\Gamma$)
a generalized Matsushima formula describing the so-called 
cuspidal cohomology in terms of the multiplicities of 
cohomological unitary representations of $G$ in 
$L^2_{\rm cusp} (\Gamma \backslash G)$, a well-known result of 
Borel \cite{BorelStable} predating \cite{FrankeSchwermer}. 
The closer study of the non-cuspidal part is the object of the 
theory of Eisenstein cohomology initiated by G. Harder 
(cf. [loc. cit.], \cite{Harder}), and there are fairly complete 
results in many cases.

While there is therefore a complete correspondence between cohomology 
and representation theory for the cuspidal part, 
generalizing the classical Eichler-Shimura homomorphism 
for the case $G = {\rm SL} (2,\mathbb{R})$ to all $G$, 
the behavior of the dimensions of the cohomology spaces is 
understood only under the hypothesis that $G$ has a compact Cartan subgroup. 
In this case, one can compute these dimensions using 
Euler-Poincar\'{e} characteristics \cite{SerreCoh} and the trace 
formula \cite{Arthur}. If the highest weight of the representation 
$E$ is regular, Arthur obtained in [loc. cit.] an 
explicit formula for the dimension of the cohomology, 
which in this case is accounted for by the packet of discrete series 
representations with the same infinitesimal character as the 
dual of $E$. In particular, the leading term in his formula is 
a constant multiple of the dimension of $E$. The question of computing the 
cohomology for non-regular highest weights and of separating the individual 
representations in the packet is connected to the problem of stabilization and to Arthur's conjectures. 
In the easiest case, namely for lattices $\Gamma$ in 
$G = {\rm SL} (2,\mathbb{R})$, it is well-known that the dimension of the cohomology 
with coefficients in the symmetric power representations of 
$G$ can be explicitly computed in terms of the basic invariants of 
$\Gamma$ (i.e. the covolume and the orders of the elliptic elements), and in fact such a dimension formula follows without 
difficulties from the description of the group-theoretic structure 
of $\Gamma$ or from the Riemann-Roch theorem (see formula (\ref{zet}) below).

The situation is different if $G$ has no compact Cartan subgroup
and there are therefore no discrete series representations. 
No explicit dimension formulas are known in this case. 
In this paper, we consider the simplest case of this type, namely 
lattices in the Lie group $G = {\rm SL} (2,\mathbb{C})$. 
Although the structure of this Lie group is very simple, 
the study of the cohomology of lattices in $G$ presents a number of 
deep problems. The irreducible finite-dimensional representations 
of $G$ are given by the tensor products
\begin{equation}\label{modu1}
E_{n,m} = \Sym^n\otimes \, \overline{\Sym}{}^m\qquad 
(n,\, m\in \Z,\ n,\, m \ge 0).
\end{equation}
Here $\Sym^n$ stands for the $n$-th symmetric power of the 
standard two-dimensional representation of $G$ and 
$\overline{\Sym}{}^m$ for its complex conjugate. 
For a lattice $\Gamma$ in $G = {\rm SL} (2,\mathbb{C})$ 
we consider therefore the finite-dimensional 
cohomology spaces
$$H^i(\Gamma,E_{n,m}).$$
The main problem studied in this paper is the behavior of the 
dimension of these spaces as a function of 
$n$ and $m$ for a fixed $\Gamma$. Another problem which we consider 
is the behavior of the dimensions 
when $\Gamma$ ranges over the subgroups of finite index in a 
lattice $\Gamma_0$. 
Since the virtual cohomological dimension of $\Gamma$ is three in the 
cocompact case, and two otherwise, 
the only dimensions with interesting cohomology 
are $i = 1$ and $i = 2$.

We now define the subspaces of cuspidal cohomology 
classes and give their description in terms of automorphic forms. 
As a consequence, it will turn out that we only have to consider the case $n=m$ and may in addition restrict to the first cohomology.
Consider the set of all proper parabolic subgroups 
$P = MU$ of ${\rm SL} (2,\mathbb{C})$ with the property that 
$\Gamma \cap U$ is a lattice in $U$. Here $U$ is the unipotent radical of 
$P$ and $M$ is a Levi subgroup of $P$.
Let $\cal C$ be a system of representatives for the 
finitely many classes of such parabolics under 
$\Gamma$-conjugation and consider the direct sum of restriction maps
\begin{equation}\label{cusp10}
H^i (\Gamma, E_{n,m}) \longrightarrow 
U^i (\Gamma, E_{n,m}) = \bigoplus_{P \in {\cal C}} 
H^i (\Gamma \cap P, E_{n,m}).
\end{equation}
The kernel of this map is called the cuspidal cohomology of 
$\Gamma$ and denoted by 
\begin{equation}
H^i_{\rm cusp}(\Gamma,E_{n,m}) \subseteq H^i(\Gamma,E_{n,m}).
\end{equation}
If $\Gamma$ is cocompact, the set $\cal C$ is empty and 
we have $H^i_{\rm cusp}(\Gamma,E_{n,m}) = H^i(\Gamma,E_{n,m})$. 

We can also describe this construction geometrically. 
The group of orientation preserving isometries of 
three dimensional hyperbolic space $X = \mathbb{H}^3$ can be 
identified with $\PSL(2,\C) = {\rm SL} (2,\mathbb{C}) / \{ \pm 1 \}$, 
and every 
lattice $\Gamma$ of  ${\rm SL} (2,\mathbb{C})$
gives rise to a quotient orbifold 
$\Gamma \backslash X$. If this orbifold is not compact, it can be 
compactified by adding a boundary $\partial(\Gamma \backslash \mathbb{H}^3)$, 
which consists 
of finitely many disjoint two-dimensional
tori or spheres. The inclusion 
$$\Gamma \backslash \mathbb{H}^3 \hookrightarrow 
\widehat{\Gamma \backslash \mathbb{H}^3}
:=\Gamma \backslash \mathbb{H}^3\cup \partial(\Gamma \backslash \mathbb{H}^3)$$
is a homotopy equivalence.
The cohomology of $\Gamma$ with coefficients in 
$E_{n,m}$ can be 
computed as the cohomology of a sheaf $\hat E_{n,m}$
on the compactified orbifold  
$\widehat{\Gamma \backslash \mathbb{H}^3}$. 
The restriction map 
\begin{equation}\label{cusp11}
H^i(\Gamma,E_{n,m})\cong
H^i(\widehat{\Gamma\backslash \mathbb{H}^3},\hat E_{n,m})\to 
H^i(\partial(\Gamma\backslash \mathbb{H}^3),\hat E_{n,m})
\end{equation}
coincides with the restriction map (\ref{cusp10}) (see \cite{GS} 
for a more detailed account).
The spaces $H^1_{\rm cusp}$ and $H^2_{\rm cusp}$ are then dual to 
each other under Poincar\'{e} duality (cf. \cite[Chapter I, \S 7]{BW}, \cite{GS}) 
and we will therefore restrict 
to the case $i=1$ in the following. 
Furthermore, by a result of 
Serre \cite[Th. 8]{SerreSL2}, for $i = 1$ the dimension of the 
image of the map in (\ref{cusp10}) is one half of the dimension of the 
target space $U^1 (\Gamma,E_{n,m})$. It is not difficult to calculate 
the latter dimension explicitly. For example, if 
$\Gamma \cap P \subseteq \pm U$ 
for all $P \in {\cal C}$ and $n+m$ is even, the dimension 
of the image is equal to the number of cusps (the number of elements of 
$\cal C$) by [loc. cit., Cor. 1]. For $i = 2$ one sees 
immediately from the long exact cohomology sequence for the pair 
$(\widehat{\Gamma \backslash \mathbb{H}^3},\partial(\Gamma 
\backslash \mathbb{H}^3))$ that the image of (\ref{cusp10}) is the 
entire target space, except in the case $n=m=0$, where it has 
codimension one. It therefore remains to study the space $H^1_{\rm cusp}$.

The theorem of Borel mentioned above yields an isomorphism 
\begin{equation}
H^1_{\rm cusp}(\Gamma,E_{n,m}) \simeq H^1 
(\mathfrak{g}, K; L^2_{\rm cusp} (\Gamma \backslash G)^\infty \otimes E_{n,m}),
\end{equation}
where the superscript $\infty$ denotes the 
subspace of smooth vectors. Since the space of cuspidal functions 
decomposes discretely as 
a representation of $G$, we can also write
\begin{equation} \label{EichlerShimura}
H^1_{\rm cusp}(\Gamma,E_{n,m}) \simeq \bigoplus_{\pi \in \hat{G}} 
{\rm Hom} (\pi, L^2_{\rm cusp} (\Gamma \backslash G)) \otimes 
H^1 (\mathfrak{g}, K; H_\pi^\infty \otimes E_{n,m}),
\end{equation}
where $\hat{G}$ denotes the unitary dual of $G$.
From the computation of the $(\mathfrak{g},K)$-cohomology of admissible irreducible representations of $G$ \cite[Chapter II]{BW}, we can deduce the following vanishing theorem:
\begin{equation}
H^1_{\rm cusp}(\Gamma,E_{n,m})=\{0\}, \quad {\rm for}\  n \neq m.
\end{equation}
In the case $n=m$ the module 
\begin{equation}
E_n: =E_{n,n}=\Sym^n\otimes \, \overline{\Sym}{}^n
\end{equation} 
is 
self-dual in the terminology of \cite{BW}, and
the dimension of $H^1_{\rm cusp} (\Gamma, E_{n})$ is equal to 
the multiplicity of the principal series representation $\pi_{2n+2,0}$ 
(the representation unitarily induced from the character 
$z \mapsto (z/|z|)^{2n+2}$ of the maximal torus $T \simeq \C^\times$) in the
space $L^2_{\rm cusp} (\Gamma \backslash G)$. We will 
use this connection extensively in the following. We would like to stress 
that it does not yield 
an explicit dimension formula. 

In the following we study the behavior of the dimension of 
$H^1(\Gamma,E_{n})$ both theoretically and numerically.
We focus primarily on the following problems.
\begin{itemize}
\item[A:] How does the dimension of $H^1(\Gamma,E_{n})$ 
behave when $\G$ is fixed and $n$ grows?
\item[B:] Are there formulas for the dimension of $H^1(\Gamma,E_{n})$
in terms of $n$ at least for some groups $\Gamma$? 
Are there formulas for the dimension valid for all $n \ge n_0$, 
where $n_0$ is allowed to depend on $\Gamma$?
\item[C:] Are there lattices $\G$ such that 
$H^1_{\rm cusp}(\Gamma,E_{n})=0$ or $H^1 (\Gamma,E_{n})=0$ for all $n$ (necessarily cocompact in the latter case)? 
\item[D:] How does the dimension of $H^1(\Gamma,E_{n})$ and 
$H^1_{\rm cusp} (\Gamma, E_{n})$
behave when $n$ is fixed and $\G$ ranges over the subgroups of finite index 
of a fixed lattice in $G$? 
\item[E:] How does the asymptotic behavior of $H^1(\Gamma,E_{n})$ and 
$H^1_{\rm cusp} (\Gamma, E_{n})$ for $n \to \infty$ change when
$\G$ ranges over the subgroups of finite index 
of a fixed lattice in $G$? 
\end{itemize}
While we do not know of previous work on problems A-C and E, problem 
D has been studied qualitatively in the context of limit multiplicities \cite{dW,Sa}
and of the conjecture of Waldhausen and Thurston (cf. \cite{DunfieldThurston}) in three-manifold topology. 
The results and computations described below are of very preliminary nature,
but we hope to provide at least some evidence on what might be true. In Section \ref{resu} below we summarize our theoretical results on upper and lower bounds for the dimension of $H^1 (\Gamma, E_n)$. In 
Section \ref{quest} we shall formulate more specific questions about 
the behavior of the these dimensions and discuss 
the numerical evidence accumulated in the later sections. 

\subsection{Theoretical results}\label{resu}

We now describe our theoretical results. 
Before we proceed, let us briefly comment on the situation for 
Fuchsian groups, i.e. discrete subgroups $\Gamma$ of 
$\SL(2,\mathbb{R})$ of finite covolume. Let $g$ be the genus of $\Gamma$, $k$ the number of cusps and $r_1$, \ldots, $r_s$ be the orders of the elliptic elements in the image of $\Gamma$ in $\PSL(2,\mathbb{R})$, considered up to conjugacy. For an integer $n$ and a positive integer $r$ let $\mu$ be the remainder of $n$ after division by $2r$ and set
\[
d (n,r) = 
\begin{cases} 1 - \frac{\mu+1}{r}, & \text{$\mu$ even}, \\
                - \frac{\mu+1}{r}, & \text{$0 \le \mu < r$ odd}, \\
              2 - \frac{\mu+1}{r}, & \text{$r \le \mu < 2r$ odd}. 
    \end{cases}
\]
Then a consideration of the group-theoretical structure of $\Gamma$ shows that
\begin{equation} \label{zet}
\dim H^1 (\Gamma, \Sym^n) = \left ( 2g - 2 + k + \sum_{i=1}^s \left( 1 - \frac{1}{r_i} \right) \right) (n+1) - \sum_{i=1}^s d (n, r_i)
\end{equation}
for all $n > 0$,
where for $-1 \in \Gamma$ one has in addition to assume $n$ to be even (the cohomology spaces vanish for odd $n$ in this case). So, the dimension of the cohomology is given by simple linear functions on congruence classes. Note also that the coefficient of $n+1$ in (\ref{zet}) is equal to ${\rm vol} (\Gamma \backslash \mathbb{H}^2) / 2 \pi$.

Our first theoretical result on lattices in $\SL (2, \C)$ concerns a general upper bound for the dimension of the cohomology.
We prove 
\begin{theorem}\label{disctheo}
Let $\Gamma \subseteq \SL(2,\C)$ be a discrete subgroup of finite covolume. Then 
$${\rm dim} \ H^1(\Gamma,E_n)=O\left(\frac{n^2}{\log n}\right)$$
as $n \to \infty$. 
\end{theorem}
This result is obtained by an application of the Selberg trace formula in 
Section \ref{dim}. 
Note that the module $E_n$ has dimension $(n+1)^2$.  
Since a group $\Gamma$ as above is  
finitely presented,
${\rm dim} \ H^1(\Gamma,E_n)=O({n^2})$ is the trivial upper bound (cf. Lemma 
\ref{triviesti} below). Non-trivial lower bounds are not known for general lattices
$\Gamma \subseteq \SL(2,\C)$. In fact, our examples in Sections \ref{hell} and \ref{tetra} 
indicate that there are probably none.

Lattices in ${\rm SL} (2,\mathbb{C})$ 
can be classified into arithmetic and 
non-arithmetic ones. The arithmetic lattices arise from 
quaternion algebras over number fields with precisely one complex place, 
and are intimately connected to number theory. 
The primary examples are the Bianchi groups ${\rm SL} (2, \O_K)$, 
where $\O_K$ is the ring of integers of an 
imaginary quadratic field $K$. It is well known that every 
non-cocompact arithmetic lattice in ${\rm SL} (2,\mathbb{C})$
 is commensurable to a Bianchi group. 
For arithmetic groups one can use Langlands functoriality 
(base change, automorphic induction, and the Jacquet-Langlands 
correspondence) to obtain lower bounds on the cohomology in certain cases
(cf. \cite{Clozel,LabesseSchwermer,Rajan}). Unfortunately, because base change is only available for solvable extensions of number fields, the results are not complete (cf. \cite[Section 6]{LackenbyLongReid}). 
The structure of the non-arithmetic lattices is even less well understood. 
The explicit examples considered in this paper are on the one 
side Bianchi groups for certain $K$ of small discriminant, and 
on the other side certain geometrically constructed lattices and series of lattices, 
almost all of which are non-arithmetic.

A related (but in general not equivalent) method for obtaining lower bounds is based on studying the action of the complex conjugation automorphism $c$ of 
${\rm SL} (2,\mathbb{C})$ on the cohomology, if the lattice $\Gamma$ is invariant under $c$ 
(as the Bianchi groups are, for example).
In this case, one can use the Lefschetz fixed point formula to compute the 
trace of this involution acting on $H^1_{\rm cusp} (\Gamma, E_n)$, 
and thereby obtain a lower bound for the dimension of this space. 
For the case of the Bianchi groups with trivial coefficients this approach 
was carried out in \cite{Kraemer,Rohlfs}. Here the results turn out to be equivalent to the ones given by theory of base change. While this method is not restricted to arithmetic groups, on the other hand it does not cover all lower bounds obtainable by Langlands functoriality for arithmetic groups. 

We work out the 
consequences of base change and automorphic induction 
(CM automorphic forms) for the cohomology of the Bianchi groups in Section \ref{change} below.
The base change construction detailed there associates to holomorphic automorphic forms for certain congruence subgroups of $\SL (2,\Z)$ elements of  
$H^1(\SL(2,\O_K),E_n)$. Let us write 
$$H^1_{\rm bc}(\SL(2,\O_K),E_n) \subseteq H^1_{\rm cusp} (\SL(2,\O_K),E_n)$$
for the corresponding subspace.
Our main result here is a precise formula for the dimension 
of $H^1_{\rm bc}(\SL(2,\O_K),E_n)$ in terms of the prime factorization of the discriminant 
of $K$. To state it, we need to introduce functions $\varepsilon_k$ and $\mu_k$ defined for all integers $k$ which depend only on the residue class of $k$ modulo $4$ and $3$, respectively, and a function $\nu_{K,k}$ depending on $K$ and at most on $k$ modulo $2$ or $3$ (see Sections \ref{changegeneral} and \ref{dimform} below for the precise definitions).

\begin{theorem}
Let $K$ be an imaginary quadratic extension of $\mathbb{Q}$, $\Ram$ the set of primes ramified in $K$ (the prime divisors of the discriminant of $K$), and for each $p \in \Ram$ let $\nu_p$ be the exact power of $p$ dividing the discriminant. 
Then there are non-negative constants $c_2$, $c_3$ and $c_4$ (depending on $K$) such that 
\begin{multline*}
\dim H^1_{\rm bc} ({\rm SL}(2,\O_K), E_n) =
\left( \frac{1}{24} \prod_{p \in \Ram} (p^{\nu_p}+1) + c_2 (-1)^{n+1} \right) 
(n+1)\\ 
 - \nu_{K,n} \frac{h_K}{2} - 2^{|\Ram|-2}
+ c_4 \varepsilon_{n+2} + c_3 \mu_{n+2} + \delta_{n,0}
\end{multline*}
for all $n \ge 0$,
where $h_K$ is the class number of $K$ and $\delta_{n,0}$ denotes the Kronecker delta symbol.
\end{theorem}

Note that for every $K$ the dimension of $H^1_{\rm bc}$ is given by linear functions of $n$ spread out over the congruence classes modulo $12$. If one takes the precise value of $c_2$ given in Section \ref{dimform} into account, one sees that the coefficient of $n$ in these linear functions is always positive, and that therefore the dimension of $H^1_{\rm cusp} (\SL(2,\O_K),E_n)$ grows at least linearly with $n$. Also, for a fixed $n$ the dimension grows linearly in the (absolute value of) the discriminant. This implies non-vanishing results for the 
cuspidal cohomology. The first results of this nature (in a much more 
limited situation) are in \cite{GSE1, GSE2}.
Let us describe the special cases 
$K = \Q (\sqrt{d})$, $d=-2$, $-7$, $-11$, more explicitly.

\begin{proposition}\label{liftintro} 
For all $n\ge 1$ we have:
\begin{equation}
\dim H^1_{\rm bc}(\SL(2,\O_{-2}),E_n) =
        \begin{cases} (n-1)/2, &\text{if $n\equiv 1\ (2)$}, \\
                      (n-2)/4, &\text{if $n\equiv 2\ (4)$},\\
                      (n-4)/4, &\text{if $n\equiv 0\ (4)$}, \\
        \end{cases}
\end{equation}
\begin{equation}
\dim H^1_{\rm bc}(\SL(2,\O_{-7}),E_n) =
        \begin{cases} (n-3)/3, &\text{if $n\equiv 0\ (3)$}, \\
                      (n-1)/3, &\text{if $n\equiv 1\ (3)$}, \\
                      (n-2)/3, &\text{if $n\equiv 2\ (3)$}, \\ 
\end{cases}
\end{equation}
\begin{equation}
\dim H^1_{\rm bc}(\SL(2,\O_{-11}),E_n) =
        \begin{cases} (n-1)/2, &\text{if $n\equiv 1\ (2)$}, \\
                      n/2, &\text{if $n\equiv 2\ (4)$}, \\
                      (n-2)/2, &\text{if $n\equiv 0\ (4)$}.\\
        \end{cases}
\end{equation}
\end{proposition}

A second construction of cohomology classes is via automorphic induction from Hecke characters of quadratic extensions of $K$ (in fact necessarily biquadratic extensions of $\Q$ unramified over $K$). In Section \ref{cmform} we describe the corresponding contribution $H^1_{\rm CM}$ to the cuspidal cohomology. In many cases, it is already contained in $H^1_{\rm bc}$. The precise criterion for an additional CM contribution to the cohomology is as follows.

\begin{proposition} \label{CMbc} Let $K$ be an imaginary quadratic field. There is a CM contribution to a space $H^1_{\rm cusp} (\SL (2,\O_K), E_n)$, $n \ge 0$, which is not contained in $H^1_{\rm bc}(\SL(2,\O_K),E_n)$, if and only if for some real quadratic field $L'$ such that $KL' / K$ is unramified, the narrow class number $h^+_{L'}$ is strictly bigger than the corresponding number of genera $g^+_{L'} = 2^{|\Ram(L')|-1}$, where $\Ram (L')$ denotes the set of primes ramified in $L'$.
\end{proposition}

See Section \ref{cmform} for the smallest examples of real quadratic fields $L'$ with this property. In any case, the additional CM contribution is always constant on residue classes modulo $12$, and the total contribution from base change and automorphic induction is therefore again given by linear functions on residue classes mod $12$.

The techniques of Langlands functoriality have been previously used to prove (a generalization of) the conjecture of Waldhausen and Thurston for some arithmetic groups $\Gamma$, i.e. to establish for fixed $n$ the existence of a finite index subgroup $\Delta$ of $\Gamma$ with $H^1 (\Delta, E_n) \neq 0$. In fact, one can in all cases in the literature show that there exists such a $\Delta$ with $H^1 (\Delta, E_n) \neq 0$ for all $n$. Building on the work of Labesse-Schwermer \cite{LabesseSchwermer}, Rajan considers in \cite{Rajan} arithmetic groups $\Gamma$ associated to quaternion algebras defined over fields $L$ such that the extension $L / L^{\rm tr}$, where $L^{\rm tr}$ is the maximal totally real subfield of $L$, is solvable (but not necessarily Galois). His methods easily imply the following result on Problem E.

\begin{proposition} Let $\Gamma$ be an arithmetic subgroup of ${\rm SL} (2, \mathbb{C})$ such that the field of definition $L$ of the corresponding quaternion algebra is a solvable extension of its maximal totally real subfield $L^{\rm tr}$. Then for every $c > 0$ there exists a finite index subgroup $\Delta$ of $\Gamma$ such that 
\[
\dim H^1 (\Delta, E_n) > c n
\]
for all $n \ge 0$.
\end{proposition}

If base change for ${\rm SL} (2)$ for arbitrary extensions of number fields was available, one could prove the corresponding result for all arithmetic lattices $\Gamma$. 

We now turn to a result which gives an upper bound for special cases 
of problem D.
We consider finite index subgroups of the Bianchi groups 
${\rm SL} (2,\O_K)$. For any non-zero ideal $\a$ of 
$\O_K$ we have the classical congruence subgroup
\begin{equation}
\Gamma_0 (\a) = \left\{\, \gamma = 
\left( \begin{array}{cc} a & b \\ c & d \end{array} \right) 
\in {\rm SL}(2,\O_K) 
\ \vrule\ \, c \in \a \right\}.
\end{equation}
The index of $\Gamma_0 (\a)$ in ${\rm SL}(2,\O_K)$ is 
given by the multiplicative function
\begin{equation}
\iota (\a) = {\rm N} (\a) \prod_{\mathfrak{p} \, | \, \a} 
\left( 1 + \frac{1}{{\rm N} (\mathfrak{p})} \right).
\end{equation}

The following theorem gives a bound for the
dimension of $H^1 (\Gamma \cap \Gamma_0 (\a), E_n)$ for each subgroup $\Gamma$
of finite index in ${\rm SL}_2 (\O_K)$, 
which improves the trivial bound $O (\iota (\a))$ by a logarithm. 

\begin{theorem} Let $\Gamma$ be a subgroup of finite index in 
${\rm SL}(2,\O_K)$, 
$K$ imaginary quadratic. Then for any fixed $n \ge 0$ we have
\begin{equation}
\dim H^1(\Gamma \cap \Gamma_0 (\a),E_n) = 
O \left(\frac{\iota (\a)}{\log {\rm N} (\a)}\right), 
\qquad {\rm N} (\a) \to \infty.
\end{equation}
\end{theorem}

This theorem, which again results from an application of the trace formula, should be compared to the limit multiplicity results of de George-Wallach \cite{dW}, L\"{u}ck and Savin \cite{Sa}, which imply 
\[
\lim_{i \to \infty} \frac{\dim H^1 (\Gamma_i, E_n)}{[\Gamma:\Gamma_i]} = 0
\]
for fixed $n$ and towers of normal subgroups $\Gamma_i$ of a fixed lattice $\Gamma$ such that $\bigcap_i \Gamma_i = \{ 1 \}$.

Note also that for $\a = a \mathfrak{o}_K$, $a$ a positive integer, 
we can get by base change arguments a lower bound 
of the form $C a = C {\rm N} (\a)^{1/2}$. If $\a$ and its 
conjugate are relatively prime, there is no non-trivial lower bound known.
On the other hand, it is certainly not possible to improve the trivial bound $O([{\rm SL} (2,\O_K):\Gamma])$ on the dimension of $H^1 (\Gamma,E_n)$ for finite index subgroups $\Gamma$ of ${\rm SL} (2,\O_K)$ without making any assumption on $\Gamma$. This follows easily from the fact that the Bianchi groups (and more generally all non-cocompact lattices in ${\rm SL} (2,\C)$) are large, i.~e. contain a finite index subgroup surjecting onto a non-abelian free group. This property is in fact conjectured to be true for all lattices (cf. \cite{LackenbyLongReid}).

\subsection{Experimental results and questions}\label{quest}

Here we formulate more specific versions of problems A-E from above.
We shall also discuss the numerical results accumulated in the later sections.

Let us begin by reporting on our numerical calculations. In Section \ref{Haeins} 
we describe how the cohomology space $H^1(\Gamma,E_n)$ can be effectively
computed from a presentation of $\G$ together with explicit matrices for the generators. 
We have developed computer codes for this task. 
The results of the computations are 
documented in Sections \ref{compu} and \ref{nonar}.
Consider first the case of Bianchi groups explained in Section \ref{compu}. We consider the fields $K = \Q(\sqrt{d})$ for
$d = -1, -2, -3, -5, -6, -7, -10, -11, -14, -19$. From Proposition \ref{CMbc} it follows immediately that
$H^1_{\rm CM} \subseteq H^1_{\rm bc}$ in all these cases. In our computations we had in all cases except two in fact $H^1_{\rm cusp} = H^1_{\rm bc}$. The precise 
range of the computations can be found in Section \ref{compu}. The two 
exceptions are:
\begin{proposition}
In the spaces 
\begin{equation}
H^1_{\rm cusp}(\SL(2,\O_{-7}),E_{12}),
\end{equation}
and
\begin{equation}
H^1_{\rm cusp}(\SL(2,\O_{-11}),E_{10})
\end{equation}
the subspace of cohomology classes obtained from base change has codimension two.
\end{proposition}
In both cases there is a uniquely determined two-dimensional complement invariant under the action of the Hecke algebra. We document the eigenvalues of the first Hecke operators on these subspaces in Section \ref{hectab} below.

These computations suggest the following question:
\begin{question} \label{BianchiQuestion}
\begin{itemize}
\item For a given field $K$, is
$$H^1_{\rm cusp}(\SL(2,\O_K),E_n)=H^1_{\rm bc}(\SL(2,\O_d),E_n)
+ H^1_{\rm CM} (\SL(2,\O_K),E_n)$$ 
for all but finitely many $n$? 
\item Is
$H^1_{\rm cusp}(\SL(2,\O_K),E_n)=H^1_{\rm bc}(\SL(2,\O_K),E_n)$  
for all $n$ for some $K$, for example $K = \Q(\sqrt{d})$, 
$d = -1, -2, -3, -5, -6, -10, -14, -19$?
\end{itemize}
\end{question}

In Section \ref{nonar} we consider some examples of (mostly) non-arithmetic lattices. 
All examples are compatible with an affirmative answer to the following question:
\begin{question} For a given lattice $\Gamma$ in ${\rm SL} (2, \mathbb{C})$, do there exist integers $n_0 \ge 0$, $N > 0$ depending on $\Gamma$ such that for each $n \ge n_0$ and $n$ in a fixed residue class modulo $N$ the dimension
${\rm dim} \ H^1(\Gamma,E_n)$ is given by a linear function in $n$?
\end{question}
As we have seen, an affirmative answer to the first part of Question \ref{BianchiQuestion} would imply that one might take $N=12$ for the Bianchi groups. A weaker but still unresolved question is:
\begin{question}\label{q1}
Do we have 
${\rm dim} \ H^1(\Gamma,E_n)=O(n)$ as  $n\to \infty$ for every lattice $\Gamma$
or is it possible that these dimensions grow faster than linearly in $n$? 
\end{question}

The computations in Section \ref{hell} and \ref{tetra} suggest an affirmative answer to the following question:
\begin{question}\label{q2}
Are there lattices $\G$ such that
${\rm dim} \ H^1(\Gamma,E_n)$ remains bounded as $n\to \infty$?
Is it possible that $H^1_{\rm cusp} (\Gamma,E_n)=0$ or even 
$H^1 (\Gamma, E_n) = 0$ for all $n$ ($\Gamma$ being necessarily cocompact in the latter case)?
\end{question}
In Section \ref{hell} an infinite sequence of non-arithmetic groups with one cusp is considered, which provides candidates for lattices with $H^1_{\rm cusp} (\Gamma,E_n) = 0$ for all $n$. In Section \ref{tetra} we consider a cocompact non-arithmetic lattice and its finite index subgroups of low index and obtain many candidates for lattices with $H^1 (\Gamma,E_n) = 0$ for all $n$.

Concerning Problem D, we pose the following variant of the conjecture of Waldhausen and Thurston as a question:
\begin{question} Given a lattice $\G$ and $n \ge 0$, is there a subgroup $\Delta$ of finite 
index in $\G$ such that 
$H_{\rm cusp}^1(\Delta,E_n) \neq 0$? 
More strongly, is there a subgroup $\Delta$ such that 
$H_{\rm cusp}^1(\Delta,E_n) \neq 0$ for all $n$?
\end{question}
We are able to provide an affirmative answer for all examples in Sections \ref{hell} and \ref{tetra} which we computed.

We have also made extensive computations of 
$\dim H^1 (\Gamma_0 (\mathfrak{p}), \C)$ for the 
standard congruence subgroups $\Gamma_0 (\mathfrak{p})$ of 
${\rm SL} (2, \O_{-1})$ associated to degree one 
prime ideals $\mathfrak{p}$ of $\O_{-1}$. The results are 
documented in Section \ref{CongruenceComp}.
The cohomology groups $H^1(\Gamma_0(\mathfrak{p}),\C)$ are particularly 
interesting for number theory since their non-vanishing is conjectured to be related to the 
existence of certain elliptic curves (or more generally abelian varieties) defined over $K=\Q(i)$ (cf. 
\cite{C, GHM, GM}). Also, the methods of Langlands functoriality do not provide any non-trivial lower bound for the dimension, and in fact there are many examples of prime ideals $\mathfrak{p}$ with $H^1(\Gamma_0(\mathfrak{p}),\C) = 0$. The analogy with distribution 
questions for elliptic curves (cf. \cite{Bru}) suggests:  
\begin{question}\label{q13} 
Is there a constant $C$ such that the asymptotic relation
$$\sum_{\mathfrak{p},\, N(\mathfrak{p})\le x} 
\dim H^1 (\G_0(\mathfrak{p}), \C)\sim C\frac{x^\frac{5}{6}}
{\log x}$$ 
holds as $x$ tends to infinity, where the sum is to be extended over all
degree one prime ideals $\mathfrak{p}$ of $\O_{-1}$ of norm at most $x$?
\end{question}
The computational results in Section \ref{CongruenceComp} are compatible with 
an affirmative answer to this question. But the range of our computations seems to be too
small to allow a more detailed analysis (cf. \cite{Bru}).
The behavior of the dimensions
$\dim H^1 (\G_0(\mathfrak{p}), E_n)$ seems to be quite different 
if $n\ge 1$ is fixed and $\mathfrak{p}$ varies (see Section \ref{CongruenceComp}).

Finally, we pose the following question regarding Problem E:
\begin{question} \label{EQuestion}
For a given lattice $\Gamma$, does there exist a subgroup $\Delta$ of finite index such that
\[
\liminf_{n \to \infty} \frac{\dim H^1 (\Delta,E_n)}{n} > 0?
\]
\end{question}
While the theoretical evidence summarized in Section \ref{resu} above suggests that this question has an affirmative answer for arithmetic lattices $\Gamma$, our computations for non-arithmetic groups are inconclusive. Namely, for the groups considered in Sections \ref{hell} and \ref{tetra} we were not able to find such a finite index subgroup $\Delta$, but to search through all subgroups of a given index very quickly becomes prohibitive.

\bigskip
{\it Acknowledgements:} We thank 
Elena Klimenko, J\"urgen Kl\"uners, 
Peter Sarnak, Haluk Sengun, Wilhelm Singhof,
Gabor Wiese and Saeid Zhargani 
for conversations on the subject.

\section{The Bianchi groups}\label{sec:groups} 

This section contains some notation and preliminary material concerning the 
Bianchi groups, as well as the explicit finite presentations on which our computer calculations are based. The first subsection fixes notation which we will 
use throughout this paper. The results needed from algebraic number theory 
are contained in \cite{La}. We also follow this book in our notational conventions.

\subsection{The Bianchi groups and their congruence subgroups}\label{biacong}

Let $d$ be a square-free negative integer, $K=\Q(\sqrt{d})\subset\C$ 
the corresponding imaginary quadratic number field and $\O=\O_d=\O_K$ 
its ring of integers. The ring $\O_d$ has a $\Z$-basis consisting of 
$1$ and $\w_d$, 
where
\begin{equation}\label{eq:omega}
  \w =\w_d =\begin{cases} \sqrt{d}, &\text{if $d\not\equiv 1 \mod 4$}, \\
                        \dfrac{1+\sqrt{d}}{2}, &\text{if $d\equiv 1\mod 4$}.
          \end{cases}
\end{equation}
The discriminant of the field $K$ is
\[ D =D_d = \begin{cases} d, &\text{if $d\equiv 1 \mod 4$}, \\
                   4d, &\text{if $d\equiv 2,3 \mod 4$}.
      \end{cases}
\]  
We set $\Ram=\Ram_d$ for the set 
of rational primes $p$  
ramified in $K$. The set $\Ram_d$ consists exactly of the prime divisors of
$D_d$. 

We also fix the following notation concerning subgroups of $\SL(2,\C)$ 
commensurable with the Bianchi groups $\SL(2,\O)$.  
Let $\a\subseteq\O$ be a non-zero ideal. The subgroup
\begin{equation}
\Gamma(\a) =\left\{\, \begin{pmatrix} a & b\\ c & d\end{pmatrix}\in
\SL(2,\O) \ \vrule \ \ a-1,\, b,\, c,\, d-1\in \a\,\right\} 
\subseteq \SL(2,\O)
\end{equation}
is called the full congruence subgroup of level $\a$. It clearly has 
finite index in $\SL(2,\O)$.
A subgroup $\Gamma\subseteq \SL(2,K)$ is called a congruence subgroup 
if $\Gamma\cap \SL(2,\O)$ has finite index in both $\Gamma$ and $\SL(2,\O)$,
and if $\Gamma$ contains a full congruence subgroup $\Gamma(\a)$ for 
a non-zero ideal $\a$ of $\O$.

Let $\a\subset K$ now  be a fractional ideal of $\O$, that is 
$\a$ is a non-zero finitely generated $\O$-submodule of $K$. We define
\begin{equation}\label{sla}
\SL(2,\a)=\left\{ \begin{pmatrix} a & b \\ c & d \end{pmatrix}\in \SL(2,K)
\ \vrule \ \ 
  a,d\in\O,\, c\in\a,\, b\in\a^{-1}\, \right\}.
\end{equation}
Notice that $\SL(2,\a)$ is a congruence subgroup of $\SL(2,K)$.
It is equal to the stabilizer in $\SL(2,K)$ of the $\O$-submodule
$\O\oplus \a$ of $K^2$.
We write $\PSL(2,\O)$ or $\PSL(2,\a)$ for the images of the corresponding
subgroups of $\SL(2,K)$ in $\PSL(2,\C)$. 

We write $\A_K$ for the ring of adeles of $K$ and $\A_{K,f}$ for the ring of 
finite adeles. We view $\A_{K,f}$ as the subring of  $\A_K$ consisting 
of those elements which are $0$ at the infinite place of $K$. The adele rings 
$\A_K$, $\A_{K,f}$ and their unit groups $\A_K^*$, $\A_{K,f}^*$
are equipped with their standard topologies (see \cite{La}). We also consider
the profinite completion of the ring $\O$, which we denote by $\hat\O$, 
to be embedded as a  
compact and open subring of $\A_{K,f}$ in the usual way.

Recall the standard description of the adelic coset space 
${\rm GL}(2,K) \backslash {\rm GL}(2,\A_K)$ in terms of
the coset spaces $\Gamma \backslash {\rm GL}(2,\C)$ for congruence 
subgroups $\Gamma$ of ${\rm GL}(2,K)$. 
By the strong approximation theorem, for
any compact open subgroup $\cal K$ of ${\rm GL}(2,\A_{K,f})$ 
the determinant map
identifies the space of connected components of
\begin{equation}
{\rm GL}(2,K) \backslash {\rm GL}(2,\A_K) / {\cal K}
\end{equation}
with the finite set
$\A_{K,f}^* / K^* \det ({\cal K})$, and for a set 
$S \subset {\rm GL}(2,\A_{K,f})$ with the property 
that $\det (S)$ forms a system of
representatives for $\A_{K,f}^* / K^* \det ({\cal K})$ we have
\begin{equation}\label{zerleg}
{\rm GL}(2,K) \backslash {\rm GL}(2,\A_K) / {\cal K} 
= \bigcup_{s \in S} \Gamma_s \backslash {\rm GL}(2,\C),
\end{equation}
where $\Gamma_s = {\rm GL}(2,K) \cap s {\cal K} s^{-1}$.

To obtain the special case of the groups ${\rm SL}(2,\a)$, let 
${\cal K}_0 = {\rm GL}(2,\hat{\O})$ be the standard maximal
compact subgroup of ${\rm GL}(2,\A_{K,f})$ and
for each finite index subgroup $\Delta$ of $\hat{\o}^*$ set 
\begin{equation}
{\cal K} (\Delta) = \{ g \in {\cal K}_0 \, | \, \det g \in \Delta \}.
\end{equation}
If $\Delta \cap \O^* = \{ 1 \}$,
the groups $\Gamma_s$ in (\ref{zerleg}) can be identified with the groups
$\SL(2,\a)$, where $\a$ runs over a system of representatives for the 
ideal classes of $K$ and each group appears with multiplicity
$[\hat{\o}^* : \Delta \o^*]$. 

Denote by $X (\Delta)$ the set of all characters of 
$\A_{K,f}^* / \Delta K^*$. It has evidently cardinality
$|X (\Delta)| = h_K [\hat{\o}^* : \Delta \o^*]$, where $h_K$ is the class number of $K$.

\subsection{Presentations}\label{present}

This subsection contains explicit finite presentations 
for some of the Bianchi groups. We include them here, because some 
of them have not yet appeared in print.
The presentations are taken from \cite{F,Sch,Sw}.
We use the standard notation for presentations of groups: 
$G=\langle \, g_1,\ldots ,g_n\ \vrule \ R_1,\ldots, R_l\, \rangle$
means that the group $G$ is generated by $g_1,\ldots,g_n$ and presented 
by the words $R_1,\ldots, R_l$. 

The following three matrices are in the set of generators in almost 
all cases:
\[ A=\begin{pmatrix} 1 & 1 \\ 0 & 1 \end{pmatrix},\quad
   B=\begin{pmatrix} 0 & 1 \\ -1 & 0 \end{pmatrix},\quad
   U=U_d=\begin{pmatrix} 1 & \w_d \\ 0 & 1 \end{pmatrix}.
\]
We first give the results for the cases $d=-1,-2,-3,-7,-11$, which are exactly the 
cases in which the
ring of integers ${\cal O}_d$ is euclidean.
\begin{equation}\label{bia1}
\PSL(2,{\cal O}_{-1})=
\left\langle \, A,B,U \ \ \vrule \ 
\ \begin{matrix} B^2,\, (AB)^3,\,  (BUBU^{-1})^3,\, AUA^{-1}U^{-1},\\ 
(BU^2BU^{-1})^2,\, (AUBAU^{-1}B)^2
\end{matrix}
\, \right\rangle,
\end{equation}
\begin{equation}\label{bia2}
\PSL(2,{\cal O}_{-2})=
\left\langle \, A,B,U \ \ \vrule \
B^2,\, (AB)^3,\, AUA^{-1}U^{-1},\, (BU^{-1}BU)^2
\, \right\rangle,
\end{equation}
\begin{equation}\label{bia3}
\PSL(2,{\cal O}_{-3})=
\left\langle \, A,\, B,\, U \ \ \vrule \ 
\ \begin{matrix} B^2,\, (AB)^3,\,  AUA^{-1}U^{-1},\\ 
(UBA^2U^{-2}B)^2, \, (UBAU^{-1}B)^3,\\ 
AUBAU^{-1}BA^{-1}UBA^{-1}UBAU^{-1}B
\end{matrix}
\, \right\rangle,
\end{equation}
\begin{equation}\label{bia7}
\PSL(2,{\cal O}_{-7})
=\langle \, A,\, B,\, U \ \vrule \ 
B^2,\, (BA)^3,\, AUA^{-1}U^{-1},\, (BAU^{-1}BU)^2\, \rangle,
\end{equation}
\begin{equation}\label{bia11}
\PSL(2,{\cal O}_{-11})
=\langle \, A,\, B,\, U \ \vrule \ 
B^2,\, (BA)^3,\, AUA^{-1}U^{-1},\, (BAU^{-1}BU)^3\, \rangle.
\end{equation}
Next we consider the case $d=-19$. In this case $\O_d$ is a 
non-euclidean principal ideal ring. We have 
\begin{equation}\label{bia19}
\PSL(2,{\cal O}_{-19})=
\left\langle \, A,\, B,\, U,\, C \ \ \vrule \ 
\ \begin{matrix} B^2,\, (AB)^3,\, AUA^{-1}U^{-1},\, C^3,\\
(CA^{-1})^3,\, (BC)^2,\, (BA^{-1}UCU^{-1})^2  
\end{matrix}
\, \right\rangle
\end{equation}
with the matrix  
$$C=\begin{pmatrix} 1-\w_{-19} & 2 \\ 2 & \w_{-19} \end{pmatrix}.$$
In the cases $d=-5,\, -6,\, -10$ the class number of $\O_d$ 
is equal to $2$. We give presentations of both  $\PSL(2,\O_d)$ and of 
$\PSL(2,\a)$ for a non-principal ideal $\a$.
\begin{equation}\label{bia5}
\PSL(2,{\cal O}_{-5})=
\left\langle \, A,\, B,\, U,\, C,\, D \ \ \vrule \ 
\ \begin{matrix} B^2,\, (AB)^3,\, AUA^{-1}U^{-1},\, D^2,\\
(BD)^2,\, (BUDU^{-1})^2,\, AC^{-1}A^{-1}BCB,\\
AC^{-1}A^{-1}UDU^{-1}CD
\end{matrix}
\, \right\rangle
\end{equation}
with matrices 
$$C=\begin{pmatrix} -4-\w_{-5} & -2\w_{-5} \\ 2\w_{-5} & -4+\w_{-5} 
\end{pmatrix}, \qquad
D=\begin{pmatrix} -\w_{-5} & 2 \\ 2 & \w_{-5} \end{pmatrix}.$$
\begin{equation}\label{bia5a}
\PSL(2,\a_{-5})=
\left\langle \, A,\, V,\, C,\, D \ \ \vrule \ 
\ \begin{matrix} AVA^{-1}V^{-1},\, CD C^{-1}D^{-1},\, (AC^{-1})^2,\\
(DV^{-1})^3,\, (CD^{-1}VA^{-1})^3
\end{matrix}
\, \right\rangle
\end{equation}
with the ideal $\a_{-5}=\langle\, 2,\, 1-\sqrt{-5}\, \rangle$ 
of $\O_{-5}$ and with the matrices
$$V=\begin{pmatrix} 1 & \frac{1+\sqrt{-5}}{2} \\ 0 & 1 \end{pmatrix},\qquad
C=\begin{pmatrix} 1 & 0 \\ 2 & 1 \end{pmatrix}, \qquad
D=\begin{pmatrix} 1 & 0 \\ 1-\sqrt{-5} & 1 \end{pmatrix}.$$
\begin{equation}\label{bia6}
\PSL(2,{\cal O}_{-6})=
\left\langle \, A,\, B,\, U,\, C,\, D \ \ \vrule \ \ 
\begin{matrix} B^2,\, (AB)^3,\,  AUA^{-1}U^{-1},\\ 
D^2,\, BCBC^{-1},\, (BAUDU^{-1})^3,\\ 
A^{-1}CAUDU^{-1}C^{-1}D^{-1},\, (BAD)^3
\end{matrix}
\, \right\rangle
\end{equation}
with the matrices
$$C=\begin{pmatrix} 5 & -2\w_{-6} \\ 2\w_{-6} & 5 \end{pmatrix},\qquad 
D=\begin{pmatrix} -1-\w_{-6} & 2-w_{-6} \\ 2 & 1+\w_{-6} \end{pmatrix}.$$
\begin{equation}\label{bia6a}
\PSL(2,\a_{-6})=
\left\langle \, 
\begin{matrix} A,\, V,\, C,\,\\ D,\, E \end{matrix} \ \ \vrule \ 
\ \begin{matrix} E^2,\, (CA^{-1})^2,\, (DV^{-1})^3,\, (DEV^{-1})^2,\\
(CEA^{-1})^2,\, CDC^{-1}D^{-1}, \\ 
AVA^{-1}V^{-1},\, (CDEV^{-1}A^{-1})^2
\end{matrix}
\, \right\rangle
\end{equation}
with the ideal $\a_{-6}=\langle\, 2,\, \sqrt{-6}\,\rangle$ of 
$\O_{-6}$ and with the matrices 
$$V=\begin{pmatrix} 1 & \frac{\w}{2} \\ 0 & 1 \end{pmatrix},\ 
C=\begin{pmatrix} 1 & 0 \\ 2 & 1 \end{pmatrix},\
D=\begin{pmatrix} 1 & 0 \\ -w & 1 \end{pmatrix}, \
E=\begin{pmatrix} -2 & -1-\frac{\w}{2} \\ 2-w & 2 \end{pmatrix}.$$
\begin{equation}\label{bia10}
\PSL(2,{\cal O}_{-10})=
\left\langle \, 
\begin{matrix} A,\, B,\, U,\, C,\\ 
D,\, E,\, F 
\end{matrix}
\ \ \vrule \ \ 
\begin{matrix} B^2,\, (AB)^3,\,  AUA^{-1}U^{-1},\, C^2,\, E^2,\\  
(BC)^2,\, (BE)^2,\, C^{-1}AD^{-1}BEBAD,\\
U^{-1}E^{-1}UFCF^{-1},\\ 
D^{-1}E^{-1}B^{-1}DU^{-1}DBCD^{-1}U,\\
D^{-1}B^{-1}ADC^{-1}U^{-1}EDA^{-1}BD^{-1}U,\\
U^{-1}DA^{-1}B^{-1}D^{-1}UFD^{-1}BADF^{-1}
\end{matrix}
\, \right\rangle
\end{equation}
with the matrices
$$C=\begin{pmatrix} -\w & 3 \\ 3 & \w \end{pmatrix},\quad
D=\begin{pmatrix} \w-1 & -4 \\ 3 & \w+1 \end{pmatrix},\quad
E=\begin{pmatrix} \w & 3 \\ 3 & -\w \end{pmatrix},\quad
F=\begin{pmatrix} 11 & 5\w \\ 2\w & -9 \end{pmatrix}.$$
\begin{equation}\label{bia10a}
\PSL(2,\a_{-10})=
\left\langle \, 
\begin{matrix} A,\, V,\, C,\,\\ D,\, E,\, F \end{matrix} \ \ \vrule \ 
\ \begin{matrix} E^2,\, (CA^{-1})^2,\, (FE)^2,\, (DEV^{-1})^2,\\ 
(DF^{-1}V^{-1})^3,\, CDC^{-1}D^{-1},\, AVA^{-1}V^{-1},\, \\
(FC^{-1}EA)^2,\, F^3,\, (CF^{-1}A^{-1})^3, \\
(CDF^{-1}A^{-1}V^{-1})^3
\end{matrix}
\, \right\rangle
\end{equation}
with the ideal $\a_{-10}=\langle\, 2,\, \sqrt{-10}\,\rangle $ 
of $\O_{-10}$ and with the matrices
$$C=\begin{pmatrix} 1 & 0 \\ 2 & 1 \end{pmatrix},\,
D=\begin{pmatrix} 1 & 0 \\ -\w & 1 \end{pmatrix},\,
E=\begin{pmatrix} -2 & -\frac{\w}{2} \\ -\w & 2 \end{pmatrix},\,
F=\begin{pmatrix} -3 & -1-\frac{\w}{2} \\ 2-\w & 2 \end{pmatrix},\, 
V=\begin{pmatrix} 1 & \frac{\w}{2} \\ 0 & 1 \end{pmatrix}.$$
The ideal class group of $\O_{-14}$ is cyclic of order 4.
The ideal $\a_{-14}=\langle\, 3,\, 1+\sqrt{-14}\, \rangle$ is not a square 
in the ideal class group. 
\begin{equation}\label{bia14}
\PSL(2,{\cal O}_{-14})=
\left\langle \, 
\begin{matrix} A,\, B,\, \\
U,\, C,\\ 
D,\, E,\, \\
F 
\end{matrix}
\ \ \vrule \ \ 
\begin{matrix} B^2,\, (AB)^3,\, 
(A^{-1}C^{-1}BDBAD^{-1}C)^2,\\  
AUA^{-1}U^{-1},\, (A^{-1}CD^{-1}ABDBC^{-1})^2,\\ 
D^{-1}CE^{-1}A^{-3}DC^{-1}A^3E,\\
CB^{-1}C^{-1}FC^{-1}BCF^{-1},\\ 
C^{-1}DA^{-1}B^{-1}D^{-1}B^{-1}CA-\\
E^{-1}A^{-2}CBD^{-1}BA^{-1}DC^{-1}A^3E, \\
ACB^{-1}D^{-1}B^{-1}A^{-1}DC^{-1}-\\
AFA^{-1}C^{-1}BDBAD^{-1}CA^{-1}F^{-1}
\end{matrix}
\, \right\rangle
\end{equation}
with the matrices
$$C=\begin{pmatrix} \w & -5 \\ 3 & \w \end{pmatrix},
D=\begin{pmatrix} 4 & 1+\w \\ 1-w & 4 \end{pmatrix},
E=\begin{pmatrix} -5+4\w & -23 \\ 4-\w & 7+\w \end{pmatrix},
F=\begin{pmatrix} 13 & 6\w \\ -2\w & 13 \end{pmatrix}.$$
\begin{equation}\label{bia14a}
\PSL(2,\a_{-14})=
\left\langle \, 
\begin{matrix} A,\, U,\, \\
C,\, D,\\ 
E,\, F,\, \\
G 
\end{matrix}
\ \ \vrule \ \ 
\begin{matrix} 
G^2,\, CDC^{-1}D^{-1},\, AUA^{-1}U^{-1},\, (CA^{-1})^3,\\  
(DGU^{-1})^2,\, F^{-1}AE^{-1}A^{-1}UFEU^{-1},\\ 
(CGE^{-1}A^{-1}UGU^{-1}AEA^{-1})^3,\\
(AEU^{-1}DGE^{-1}A^{-1}UGD^{-1})^2,\\ 
DC^{-1}GU^{-1}AEGD^{-1}UE^{-1}F^{-1}-\\
CGE^{-1}A^{-1}UGU^{-1}AEA^{-1}F
\end{matrix}
\, \right\rangle
\end{equation}
with the matrices
$$U=\begin{pmatrix} 1 & \frac{1-\w}{3} \\ 0 & 1 \end{pmatrix},\ 
C=\begin{pmatrix} 1 & 0 \\ 3 & 1 \end{pmatrix}, \
D=\begin{pmatrix} 1 & 0 \\ 1+\w & 1 \end{pmatrix}, \
E=\begin{pmatrix} -3-\w & -4 \\ 6 & 3-\w \end{pmatrix},$$
$$F=\begin{pmatrix} -3+\w & -3 \\ 2+2\w & -3+\w \end{pmatrix},\quad 
G=\begin{pmatrix} -2 & \frac{\w-1}{3} \\ 1+\w & 2 \end{pmatrix}.$$

\section{Group cohomology}\label{cohomol}

In this section we report basic definitions from the cohomology of groups.
Section \ref{Haeins} reports a method to compute the first cohomology 
group for finitely presented groups. Our basic reference here is \cite{B}. 

Let $\G$ be a group and $M$ a $R\G$-module for a commutative ring $R$. 
A derivation from $\G$ to $M$ is a map $f:\G\to M$ which satisfies
\begin{equation}\label{eq:der} 
 f(gh)=g\cdot f(h)+f(g)
\end{equation}
for all $g,h\in \G$. For $m\in M$ the map 
\begin{equation}\label{eq:ider} 
f_m: \G\to M,\qquad f_m(g) =g\cdot m-m, 
\end{equation}
is a derivation and is called the inner derivation corresponding to $m$.
We write $\Der(\G,M)$ for the space of all derivations and  
$\IDer(\Gamma,M)$ for its subspace consisting of inner derivations. 
If $H^1(\G,M)$ is the first cohomology group of $\G$ with coefficients 
in $M$,
we have
\begin{equation}\label{dere}
H^1(\G,M)=\Der(\G,M)/\IDer(\G,M). 
\end{equation}

Here we are interested in the case when $\G\subseteq \SL(2,L)\subseteq \SL(2,\C)$ where 
$L\subset\C$ is a number field. The modules we consider are derived from the 
symmetric powers of the standard representation of $\SL(2,\C)$. So, let
$V$ be a two-dimensional $L_1$-vector space with basis $x,\, y$ where 
$L_1$ is a field between $L$ and $\C$ invariant under complex conjugation.
Let
$n$ be a non-negative integer. 
The symmetric power $\Sym^n(L_1)$ has the $L_1$-basis
$x^{n-i}y^i$, $0\le i\le n$.
The action of $g\in\SL(2,L_1)$ is given by  
\begin{equation}\label{actio} 
g\cdot x^{n-i}y^i=\begin{pmatrix} a & b \\ c & d \end{pmatrix} 
\cdot x^{n-i}y^i= (ax+cy)^{n-i}(bx+dy)^i,
\qquad g=\begin{pmatrix} a & b \\ c & d \end{pmatrix}.
\end{equation}
The module $\overline{\Sym}{}^n(L_1)$ is equal to 
${\Sym}{}^n(L_1)$ as an $L_1$-vector space, and the action is
given by replacing $g$ in (\ref{actio}) by its  
complex conjugate.

We often use the follwing simple facts from group cohomology without
further notice. First of all the spaces 
$H^1(\G,{\Sym}{}^n(L_1)\otimes \overline{\Sym}{}^m(L_1)) \otimes \C$
and 
$H^1(\G,{\Sym}{}^n(\C)\otimes \overline{\Sym}{}^m(\C)))$ are isomorphic
for all $n,\, m\ge 0$. Secondly, if $m+n$ is even, the action of $\G$ 
on ${\Sym}{}^n(L_1)\otimes \overline{\Sym}{}^m(L_1)$ factors through 
an action of the image $\tilde\G$ of $\G$ in $\PSL(2,\C)$ and 
$H^1(\G,{\Sym}{}^n(L_1)\otimes \overline{\Sym}{}^m(L_1))$ 
is isomorphic to
$H^1(\tilde\G,{\Sym}{}^n(L_1)\otimes \overline{\Sym}{}^m(L_1))$.

\subsection{$H^1(\G,M)$ for finitely presented groups}\label{Haeins}

Here we explain how information about $H^1(\G,M)$ can be computed from 
equation (\ref{dere}). We assume here that $R$ is an euclidean ring and $M$ is
a free $R$-module of finite rank in which a basis has been chosen. 
Let $\G$ be a finitely presented group given explicitly in the form
$$\G=\langle \, g_1,\dots,g_s\ \vrule  \ R_1,\ldots, R_t\, \rangle.$$ 
Here we consider the 
relations $R_1,\dots,R_t$ to be explicitly given words in the 
generators $g_1,\dots,g_s$ of $\G$ and their inverses. Assume also that
the matrices for the action of $g_1,\dots,g_s$ on $M$ are explicitly given.
Consider now the $R$-linear map
$$\Phi: {\rm Der}(\G,M)\to M^s,\qquad \Phi(f) =(f(g_1),\ldots,f(g_s)).$$
The image of $\Phi$ lies in the kernel of the linear map 
$\Lambda: M^s\to M^t$, which is obtained by formally expanding 
the image of each of the relators $R_1,\ldots, R_t$ under a derivation 
$f:\G\to M$ in terms of the values $f(g_1),\ldots,f(g_s)$. It is easily seen that
$\Phi$ maps ${\rm Der} (\G,M)$ isomorphically to the kernel ${\rm ker}(\Lambda)$ of $\Lambda$.
Since $M^s$ is a free $R$-module, a basis for the free module 
${\rm ker}(\Lambda)$ can be computed.
Consider now the linear map
$$\mu: M\to {\rm ker}(\Lambda),\qquad \mu(m) =((g_1-1)m,\ldots,(g_s-1)m).$$
The image of $\mu$ may then be described as the linear span of the 
images of the basis elements of $M$. If we express these 
in terms of the previously computed basis of ${\rm ker}(\Lambda)$, we see that the
effective version of the elementary divisor theorem can be used to 
compute the structure of    
\begin{equation}\label{dereco}
H^1(\G,M)=\Der(\G,M)/\IDer(\G,M) = {\rm ker}(\Lambda)/{\rm Im}(\mu).
\end{equation}
If $R$ is a field, the dimension of $H^1(\G,M)$ can be computed 
by this method.

Apart from being important for the computation of 
cohomology spaces, (\ref{dereco}) leads to the following (trivial) estimate.
\begin{lemma}\label{triviesti}
Let $\G$ be a group 
generated by $s$ elements and $M$ a finite dimensional 
$R\G$-module for some field $R$. Then
$\dim H^1(\G,M)\le s\dim M.$
\end{lemma}
A typical problem encountered in our computations is that the module $M$
can be a vector space of big dimension (up to around $50000$) over 
an algebraic number field, and that the direct computation of   
the dimension of $H^1(\G,M)$ 
from (\ref{dereco}) is not feasible. 
All discrete subgroups $\G \subseteq \SL(2,\C)$ considered in this paper 
have the property that they are contained in $\SL (2,R)$ for a finitely generated
ring $R$ inside an algebraic number field. Let $\O_\G$ be a ring containing $R$ and its complex conjugate. 
Suppose $p$ is a prime and $\O_\G\to \F_p$ is a surjective 
ring homomorphism. Then $E_n(\F_p)$ inherits
the structure of a $\G$-module. By the usual universal coefficient theorem we have
$$\dim_{\F_p}H^1(\G,E_n(\F_p))\ge \dim_{\C}\, H^1(\G,E_n).$$
A standard argument using Tchebotarev's density theorem shows that
$\dim_{\C}\, H^1(\G,E_n)$ is equal to the minimum of the dimensions
$\dim_{\F_p}H^1(\G,E_n(\F_p))$, where $p$ ranges over all primes with 
the above compatibility property. For all real numbers $x$ we define  
\begin{equation}
\dim_{\le x}\, H^1(\G,E_n) = \inf_{p \le x}\  \{\, \dim_{\F_p}H^1(\G,E_n(\F_p)) \,\},  
\end{equation}
where $p$ ranges over all primes with $p\le x$ which admit 
a surjective ring homomorphism $\O_\G\to \F_p$. 
The numbers $\dim_{\le x}\, H^1(\G,E_n)$ are much cheaper to
compute than the actual dimensions $\dim_{\C}\, H^1(\G,E_n)$.
Of course, in the computations below 
we hope to have chosen the bound $x$ to be large enough to 
capture $\dim_{\C}\, H^1(\G,E_n)$. Also, if a lower bound for this dimension is known beforehand, we can by this method verify that the actual dimension is equal to the bound.

\subsection{Hecke operators}\label{hecke} 

In this section we introduce the Hecke operators on cohomology spaces in a way suitable for explicit computations. We chose
a treatment similar to \cite[Section 8.5]{Sh1}, see also \cite{GHM}.

If $H$ is a subgroup of a group $\G$, 
and $M$ is a $\G$-module, the inclusion $H\hookrightarrow \G$
induces a restriction map 
$\res^\G_H: H^*(\G,M)\longrightarrow H^*(H,M)$.
When $[\G:H]<\infty$, there is also a map 
$\tr: H^*(H,M)\longrightarrow H^*(\G,M)$ 
in the opposite direction, called the transfer map (cf. \cite{B}).
The composition 
$\tr\circ\res^\G_H$ is multiplication by
$[\G:H]$ on $H^*(\G,M)$.

Let now $\G$ be a congruence subgroup of $\SL(2,\O)$, where $\O$ is the 
ring of integers in an imaginary quadratic number field $K$, and $M$
be one of the $\GL(2,\C)$-modules $E_{n,m}$.  
The groups $\G$ and $\delta\G\delta^{-1}$ are easily seen to be 
commensurable for every $\delta\in \GL(2,K)$. 
Define the Hecke operator $T_{\delta}: H^1(\G,M)\to H^1(\G,M)$ 
by the diagram:
\begin{equation}
 \begin{CD}
  H^1(\G,M)  @> T_{\delta} >> H^1(\G,M) \\
  @V \res VV   @AA \tr A \\
  H^1(\G\cap \delta\G\delta^{-1},M) @> \tilde{\delta} >> 
H^1(\delta^{-1}\G\delta\cap \G,M)
 \end{CD}
\end{equation}
where $\tilde\delta$ Is the isomorphism in cohomology induced by 
conjugation with $\delta$. For a non-zero element $a \in \O$ we define
\begin{equation}\label{hedef}
T_a =T_{\delta_a}\qquad {\rm with}\qquad 
\delta_a=\begin{pmatrix} 1 & 0 \\ 0 & a \end{pmatrix}.
\end{equation}
The following properties of the linear maps 
$T_{\delta}: H^1(\G,M)\to H^1(\G,M)$, $\delta\in \GL(2,K)$, are well
known (cf. \cite[Section 8.5]{Sh1}, \cite{GHM}):
\begin{itemize}
\item Each $T_\delta$ is diagonalizable.
\item The characteristic polynomial of $T_\delta$ has integral coefficients 
and its zeroes are real numbers.
\item $T_\delta$ depends only on the double coset $\G\delta\G$.
\item If $\G=\SL(2,\O)$, all operators $T_\delta$ commute with each other.
\end{itemize}

\subsection{The Eichler-Shimura isomorphism}\label{ESH}

In this subsection we briefly recall the generalized Eichler-Shimura isomorphism sketched already in the introduction,
which will give us the possibility of using results from the theory of automorphic forms in our study of cohomology spaces. See also \cite{Harder} and \cite[Th\'{e}or\`{e}me 3.2]{Urban} for the case of congruence subgroups of $\GL (2, K)$, $K$ imaginary quadratic.

From \cite[Chapter II]{BW} we know that for any integer $n \ge 0$, and any unitary representation $\pi$ of $G = \SL (2, \C)$, the $(\mathfrak{g},K)$-cohomology space $H^1 (\mathfrak{g},K; H_\pi^\infty \otimes E_n)$ is non-trivial if and only if $\pi$ is the principal series representation $\pi_{2n+2,0}$ (the representation unitarily induced from the character 
$z \mapsto (z/|z|)^{2n+2}$ of the maximal torus $T \simeq \C^\times$, cf. Section \ref{dim}), and one-dimensional in this case. Therefore, we can deduce from (\ref{EichlerShimura}) the more explicit isomorphism
\begin{equation}
\Hom (\pi_{2n+2,0}, L^2_{\rm cusp} (\Gamma \backslash {\rm SL}(2,\C)) 
\simeq H_{\rm cusp}^1 (\Gamma, E_n)
\end{equation}
for any lattice $\Gamma$ of $G$.

For use in Section \ref{change}, we quickly rewrite this isomorphism in a form involving $\GL (2, \C)$.
Define a unitary character of $\C^*$ by $\chi_\infty (x) = x/|x|$, and 
for each integer $n \ge 0$ 
consider the principal series representation 
$\rho^n_\infty = {\rm PS} (\chi_{\infty}^{n+1},
\chi_{\infty}^{-n-1})$ of
${\rm GL}(2,\C)$. Let $Z_\infty \subset {\rm GL}(2,\C)$ be the center of
${\rm GL}(2,\C)$. Then we have an isomorphism
\begin{equation}
\Hom (\rho_n^\infty, L^2_{\rm cusp} (\Gamma \backslash {\rm GL}(2,\C) / Z_\infty)) 
\simeq H_{\rm cusp}^1 (\Gamma, E_n).
\end{equation}

\section{Base change}\label{change}

This section contains our results on the construction of 
cohomology classes for the Bianchi groups by base change from classical modular forms for congruence subgroups of $\SL(2,\Z)$ and automorphic induction from Hecke characters of quadratic extensions. In particular, we derive explicit dimension formulas for the corresponding subspaces of the cohomology.
For this we fix an imaginary quadratic number field $K=\Q(\sqrt{d})$ and use 
the notation of Section \ref{biacong}.

\subsection{General results on the base change construction} \label{changegeneral}

Here we present the consequences of the theory of base change and automorphic induction
for the
cohomology of the groups $\SL (2,\a)$. 
We give a precise description of the base change process and the relevant spaces of holomorphic elliptic modular forms.
The notation and 
concepts from the theory of automorphic forms are taken from \cite{BW}.
For a quadratic extension $L$ of $\Q$ 
denote by $\omega_L$ the associated quadratic character
of $\A_\Q^* / \Q^*$. 

Let $\Autf_K$ be the set of all cuspidal 
automorphic representations of ${\rm GL}(2,\A_K)$.
See \cite{Langlands} for information on the 
base change map $\pi \mapsto \pi_K$ from
${\rm GL}(2,\A_\Q)$ to ${\rm GL}(2,\A_K)$.
We shall be interested in the following subset of $\Autf_K$.
\begin{definition}
The set $\Autf_K^{bc}$ of (twisted) base change 
representations is the set of all 
$\Pi \in \Autf_K$ such that $\Pi \simeq \pi_K \otimes \chi$ 
for an automorphic representation
$\pi$ of ${\rm GL}(2,\A_\Q)$ and an idele class 
character $\chi$ of $K$. 
\end{definition}
Recall from Section \ref{ESH} the definition of the representations
$\rho^n_\infty$ of $\GL (2, \C)$. For an integer $n \ge 0$ and a finite 
index subgroup $\Delta$ 
of $\hat{\o}^*$ consider
\begin{equation}
\Autf^1_K (n, \Delta) = \{ \Pi \in \Autf_K \, | \, 
\Pi_\infty \simeq \rho^n_\infty, \, \Pi_f^{{\cal K} (\Delta)} \neq 0 \}
\end{equation}
and set
\begin{equation}
\Autf^1_K (n) = \bigcup_{\Delta} \Autf^1_K (n, \Delta).
\end{equation}
Furthermore, let 
$\Autf^{1,bc}_K (n, \Delta) = \Autf^1_K (n, \Delta) \cap \Autf_K^{bc}$ and 
$\Autf^{1,bc}_K (n) = \Autf^1_K (n) \cap \Autf_K^{bc}.$
Recall that each representation in $\Autf_K$ occurs with multiplicity one in 
$L^2_{\rm cusp} ({\rm GL}(2,K) \backslash {\rm GL}(2,\A_K))$. Furthermore, $\Pi \in \Autf^1_K (n)$ is equivalent to the condition that 
$\Pi_\infty \simeq \rho^n_\infty$ and that the local components $\Pi_{\mathfrak{p}}$ at the finite places $\mathfrak{p}$ are twists of unramified principal series representations by characters. Therefore, for $\Pi \in \Autf^1_K (n, \Delta)$ the space $\Pi_f^{{\cal K} (\Delta)}$ is actually one-dimensional.

If we take a subgroup $\Delta$ of $\hat{\o}^*$ with the property $\Delta \cap \o^* = \{ 1 \}$, 
we have by (\ref{zerleg}) an isomorphism
\begin{multline*}
\left( \bigoplus_\a \Hom (\rho_n^\infty, 
L^2_{\rm cusp} (SL(2,\a) \backslash {\rm GL}(2,\C) / Z_\infty)) 
\right)^{[\hat{\o}^* : \Delta \o^*]} \\
\simeq  \Hom (\rho_n^\infty, L^2_{\rm cusp}
({\rm GL}(2,K) \backslash {\rm GL}(2,\A_K) / Z_\infty {\cal K} (\Delta))),
\end{multline*}
where $\a$ ranges over a system of 
representatives for the ideal classes of $K$.
Combining the Eichler-Shimura isomorphism from Section \ref{ESH} with multiplicity one and the fact that $\dim \Pi_f^{{\cal K} (\Delta)} = 1$ for $\Pi \in \Autf^1_K (n, \Delta)$,
we obtain the relation
\begin{equation}
\sum_\a \dim H_{\rm cusp}^1 (\SL(2,\a), E_n)
=  \frac{\left| \Autf^1_K (n, \Delta) \right|}
{[\hat{\o}^* : \Delta \o^*]}.
\end{equation}
It is not difficult to obtain also 
a finer description distinguishing between the individual cohomology spaces
$\dim H_{\rm cusp}^1 (\SL(2,\a), E_n)$ 
for representatives $\a$ of different ideal classes. 
For this consider the action of
the abelian group $X (\Delta)$ on the space 
$\Hom (\rho_n^\infty, L^2_{\rm cusp}
({\rm GL}(2,K) \backslash {\rm GL}(2,\A_K) / 
Z_\infty {\cal K} (\Delta)))$
 given by letting $\xi \in X (\Delta)$ (see Section \ref{biacong})
act as multiplication of functions 
on 
$${\rm GL}(2,K) \backslash {\rm GL}(2,\A_K) / 
Z_\infty {\cal K} (\Delta)$$ 
by $\xi \circ \det$. 
Considering a basis of 
$\Hom (\rho_n^\infty, L^2_{\rm cusp}
({\rm GL}(2,K) \backslash {\rm GL}(2,\A_K) / 
Z_\infty {\cal K} (\Delta)))$ consisting
of normalized (cf. \cite[Section 5]{Urban}) eigenfunctions for the Hecke algebra
of ${\cal K} (\Delta)$ (which correspond to the representations 
in $\Autf^1_K (n, \Delta)$),
one sees that the action of 
$X (\Delta)$ induces a permutation of this basis, and therefore
the trace of the action of a non-trivial element 
$\xi \in X (\Delta)$ is equal to the number of 
elements $\Pi \in \Autf^1_K (n, \Delta)$
with $\Pi \otimes \xi \simeq \Pi$. 
It is clear that this number can only be non-zero if $\xi$ is quadratic, 
and indeed unramified quadratic, 
i.~e. necessarily of the form $\omega_L \circ {\rm N}_{K/\Q}$ 
for an imaginary quadratic field
$L \neq K$ such that $LK/K$ is unramified 
(see Proposition \ref{CMclass} below). 
Therefore, we get
\begin{multline*}
\dim H_{\rm cusp}^1 (\SL(2,\a), E_n) 
= \frac{1}{|X(\Delta)|}
\Biggl(
\left| \Autf^1_K (n, \Delta) \right| + \mbox{}  \\
\sum_{L \in {\cal L} (K)} \omega_L ({\rm N}_{K / \Q} (\a)) 
\left| \{ \Pi \in \Autf^1_K (n, \Delta) \, | \, 
\Pi \otimes \omega_L \circ {\rm N}_{K / \Q} \simeq \Pi \} \right|
\Biggr),
\end{multline*}
where ${\cal L} (K)$ denotes the set of all 
imaginary quadratic fields $L \neq K$ with $LK/K$ unramified.
Furthermore, if $A \subseteq \Autf^1_K (n, \Delta)$ is 
any subset invariant under twisting by characters in $X(\Delta)$, we can
consider inside the space $L^2_{\rm cusp} ({\rm GL}_2 (K) 
\backslash {\rm GL}(2,\A_K) / Z_\infty {\cal K} (\Delta))$ the subspace
spanned by representations in $A$ and apply the 
same arguments to see that it splits as a direct sum of spaces of functions
supported on a single connected component. 
This implies that it makes sense to speak of the 
contribution of representations in $A$ 
to each space $H_{\rm cusp}^1 (\SL (2, \a), E_n)$ 
and that the 
dimension of the corresponding subspace is given by
\begin{multline} \label{H1Gammaa}
\dim H_{{\rm cusp},A}^1 (\SL(2,\a), E_n) \\
= \frac{1}{|X(\Delta)|}
\left(
|A| + \sum_{L \in {\cal L} (K)} \omega_L ({\rm N}_{K / \Q} (\a))
\left| \{ \Pi \in A \, | \, \Pi \otimes 
\omega_L \circ {\rm N}_{K / \Q} \simeq \Pi \} \right|
\right).
\end{multline}
In particular, this dimension depends only on the genus of $\a$ and
it assumes its maximum on the principal genus. 
We are especially interested in evaluating the 
contribution of twisted base change forms to 
the cohomology, i.~e. in the case $A = \Autf^{1,bc}_K (n, \Delta)$.

\begin{definition}
For the set $A_{\rm bc}  = \Autf^{1,bc}_K (n, \Delta)$ define
\[
H^1_{\rm bc} (\SL(2,\a), E_n) := H_{{\rm cusp},A_{\rm bc}}^1 (\SL(2,\a), E_n) \subseteq H_{{\rm cusp}}^1 (\SL(2,\a), E_n).
\]
\end{definition}
Note that this definition makes sense, since the right-hand side is indeed independent of the subgroup $\Delta$ with $\Delta \cap \O^* = \{ 1 \}$.

Our first goal is to describe the set
$\Autf^{1,bc}_K (n, \Delta)$ in terms of holomorphic 
automorphic forms for ${\rm GL}(2,\A_\Q)$ 
fulfilling explicit local conditions.
We also need to distinguish the automorphic representations 
of CM type. Recall that for any quadratic extension $E/F$ of 
number fields there is a canonical map from Hecke characters of 
$E$ to automorphic representations of ${\rm GL}(2,\A_F)$ 
\cite{JacquetLanglands} which is called automorphic induction 
(notation: ${\rm AI}_{E/F}$). The map is characterized by 
${\rm AI}_{E/F} (\theta \, \chi \circ {\rm N}_{E/F}) = 
{\rm AI}_{E/F} (\theta) \otimes \chi$
for Hecke characters $\chi$ of $F$, and the 
$L$-function identity $L (s, {\rm AI}_{E/F} (\theta)) = 
L (s, \theta)$. The automorphically induced representation 
is cuspidal if and only if $\theta^\tau \neq \theta$, 
where $\tau$ is the automorphism of $E/F$, and the 
fibers of the automorphic induction map are precisely the 
orbits $\{ \theta, \theta^\tau \}$ of $\tau$. 
There are compatible local induction maps, which we also denote by ${\rm AI}$.

In the following, we fix once and for all for each $p \in \Ram$ a 
character $\theta_p$ of $K^*_p$ with
$\theta_p / \theta_p^c$ unramified quadratic. 
\begin{definition} \label{loccomp}
For each $n \ge 0$ let ${\cal A}^1_{\Q} (n)$ be the set of all 
cuspidal automorphic
representations $\pi$ of ${\rm GL}(2,\A_\Q)$ such that $\pi_\infty$ is the
holomorphic discrete series representation of weight $n+2$, 
$\pi_p$ is unramified for $p \notin \Ram$, and of one of the
following three types for $p \in \Ram$:
\begin{enumerate}
\item
unramified principal series,
\item ${\rm PS} (\alpha, \omega_{K,p} \beta)$ with $\alpha$, $\beta$ 
unramified characters of $\Q^*_p$,
\item ${\rm AI}_{K_p/\Q_p} (\theta_p) \otimes \gamma$ with an 
unramified character $\gamma$ of $\Q^*_p$.
\end{enumerate}
\end{definition}

For any (necessarily imaginary) quadratic
extension $L$ of $\Q$ let $\Autf^1_{\Q} (n; L)$ be the subset of 
$\Autf^1_{\Q} (n)$ consisting of
representations automorphically induced from $L$. 
Recall that $\pi \in \Autf^1_{\Q} (n; L)$ 
if and only if $\pi \otimes \omega_L \simeq \pi$ \cite{SL2}.
We will see that the set of possible extensions 
$L$ is precisely ${\cal L} (K)$.

The basic classification statement is the following proposition. It shows that
we obtain the representations in $\Autf^{1,bc}_K (n, \Delta)$ by base change and character twists from the 
elliptic modular forms satisfying the local conditions of Definition \ref{loccomp}. Of course, the description depends on the choice of the local characters $\theta_p$ for $p \in \Ram$.

\begin{proposition}
\begin{enumerate}
\item
If $\Pi \in \Autf^{1,bc}_K (n)$, one can find $\pi \in {\cal A}^1_{\Q} (n) 
\backslash {\cal A}^1_{\Q} (n; K)$
such that $\Pi \simeq \pi_K \otimes \chi$ for some finite 
order idele class character $\chi$ of $K$.
\item
If $\Delta \cap \o_K^\times = \{1\}$, for any 
$\pi \in {\cal A}^1_{\Q} (n) \backslash {\cal A}^1_{\Q} (n; K)$
there exists a finite order idele class
character $\chi$ of $K$ with $\pi_K \otimes \chi 
\in \Autf^1_K (n, \Delta)$. The set of all such characters $\chi$ is
a principal homogeneous space for $X (\Delta)$.
\end{enumerate}
\end{proposition}

{\em Proof.} For an automorphic representation $\pi$ of ${\rm GL}(2,\A_\Q)$ we have
$\pi_K \in \Autf^1_K (n)$ if and only if $\pi$ is not automorphically induced from $K$, $\pi_\infty$ is up to a character twist the holomorphic or anti-holomorphic discrete series representation of weight $n+2$, and each $\pi_p$ satisfies up to a character twist the local condition of Definition \ref{loccomp}. It is not difficult to deduce from this the assertions of the Proposition.

We can also classify the CM representations of interest to us as follows.

\begin{proposition} \label{CMclass}
\begin{enumerate}
\item If $\Autf^1_{\Q} (n; L)$ is nonempty, 
$L$ is an imaginary quadratic extension of $\Q$ such that
for all primes $p$ the character $\omega_{L,p}$ is either 
unramified or the product of $\omega_{K,p}$ and an
unramified character.  
\item
If for $\Pi \in \Autf^1_K (n)$ there exists a character 
$\gamma \neq 1$ with $\Pi \otimes \gamma \simeq \Pi$, 
the character $\gamma$ is necessarily of the form 
$\omega_{KL/K} = \omega_L \circ {\rm N}_{K/\Q}$ for some
quadratic extension $L / \Q$ as above.
\end{enumerate}
\end{proposition}

The set of all imaginary quadratic number fields different 
from $K$ and satisfying the conditions of Part 1 of 
Proposition \ref{CMclass} is precisely the set 
${\cal L} (K)$ of imaginary quadratic fields different from 
$K$ for which $LK /K$ is unramified. Equivalently, 
it is the set of all imaginary quadratic fields $L$ 
for which the discriminant $d_L$ is a proper 
divisor of the discriminant of $K$ and the two 
factors $d_L$ and $d_K / d_L$ are coprime.

Consider now for each $\pi \in \Autf^1_{\Q} (n)\backslash 
{\cal A}^1_{\Q} (n; K)$ the set
\[
\Autf^1_K (n, \Delta; \pi) := \{ \Pi = \pi_K \otimes 
\chi \, | \, \Pi \in \Autf^1_K (n, \Delta) \}.
\]
Clearly, the sets $\Autf^1_K (n, \Delta; \pi)$ form a partition of $\Autf^{1,bc}_K (n, \Delta)$.
Assuming $\Delta \cap \o_K^* = \{1\}$, the set $\Autf^1_K (n, \Delta; \pi)$ 
has cardinality 
$|X (\Delta)|$, if $\pi$ is not automorphically 
induced from any quadratic field $L$, and 
$|X (\Delta)|/2$, otherwise. It remains to 
count for any $\pi \in \Autf^1_{\Q} (n)$ the number of 
$\pi'$ with $\Autf^1_K (n, \Delta; \pi') = \Autf^1_K (n, \Delta; \pi)$. 
We first consider the non-CM representations.

\begin{proposition} 
For $\pi \in \Autf^1_{\Q} (n) \backslash 
\bigcup_{L \in {\cal L} (K) \cup \{ K \}} {\cal A}^1_{\Q} (n; L)$ the set 
\[
\{ \pi' \in \Autf^1_{\Q} (n) \, | \, 
\Autf^1_K (n, \Delta; \pi') = \Autf^1_K (n, \Delta; \pi) \}
\]
consists of the twists $\pi \otimes \gamma$ for all 
characters $\gamma$ such that
$\gamma_p$ is unramified for all $p$ where $\pi_p$ is unramified, 
and $\gamma_p$ is unramified or the product of $\omega_{K,p}$ 
and an unramified character at the primes
$p$ where $\pi_p$ is ramified. In particular, it has
cardinality $2^{|R(\pi)|}$, where $R(\pi) \subseteq \Ram$ 
denotes the set of all primes $p$ where $\pi_p$
is ramified.
\end{proposition}

Therefore, if we want to write the cardinality of $\Autf^{1,bc}_K (n, \Delta)$ as a sum over all representations $\pi \in \Autf^1_{\Q} (n)\backslash 
{\cal A}^1_{\Q} (n; K)$, each non-CM representation $\pi \in \Autf^1_{\Q} (n) \backslash 
\bigcup_L {\cal A}^1_{\Q} (n; L)$ has to be weighted by the factor
$|X(\Delta)| 2^{-|R(\pi)|}$.

\begin{example}
Consider the case where a single prime $p$ is ramified in $K$. 
In this case, the set
$\Autf^1_{\Q} (n)$ consists of the automorphic representations 
associated to classical modular forms
of weight $n+2$ for ${\rm SL}(2,\mathbb{Z})$, for $\Gamma_0 (p)$ 
with character $\omega_K$, or of $p$-power level with
$\pi_p \simeq {\rm AI} (\theta_p) \otimes \gamma_p$, 
$\gamma_p$ unramified. The CM forms for $K$ have to be
omitted. In this case, there are no other fields $L$ to be considered. 
To obtain the
dimension of $H^1_{\rm bc} (\SL (2, \O), E_n)$, 
the dimension of the corresponding spaces of modular forms
has to be weighted by a factor $1 / 2$ except in the ${\rm SL}(2,\mathbb{Z})$ 
case.
\end{example}

In the count for the representations in $\Autf^{1,bc}_K (n, \Delta)$,
the main term is therefore given by
$|X(\Delta)| \sum_{\pi \in \Autf^1_\Q (n)} 2^{-|R(\pi)|}$.
The contributions from CM representations have to be modified 
by omitting the representations automorphically induced from the field $K$ and
weighting the contribution of the representations induced 
from quadratic fields $L \in {\cal L} (K)$ by an additional factor $\frac{1}{2}$. The
reason for this is that for these representations there are more equivalences
$\Autf^1_K (n, \Delta; \pi') = \Autf^1_K (n, \Delta; \pi)$ than 
in the non-CM case.

To give some more details, we first explicate the local 
conditions on CM representations 
in $\Autf^1_{\Q} (n)$. Recall the definition of the local character at infinity $\chi_\infty$ in Section \ref{ESH}.

\begin{lemma}
Let $L \in {\cal L} (K)$ be an imaginary quadratic field. Write
${\rm AI}_{K_p / \Q_p} (\theta_p) = {\rm AI}_{L_p / \Q_p} 
(\theta_{p,L_p})$ with 
a character $\theta_{p,L_p}$ of $L_p^\times$ for all 
$p \in \Ram$ where $p$ is non-split in $L$
(note that this is possible). For an idele class character 
$\psi$ of $L$ with $\psi_\infty = \chi_\infty^{-(n+1)}$ and
unramified at primes not above primes in $\Ram$
we have ${\rm AI}_{L/\Q} (\psi) \in \Autf^1_{\Q} (n)$ if and 
only if the following local conditions are satisfied:
\begin{enumerate}
\item If $p \in \Ram$ splits in $L$, $\psi_p$ is either unramified 
or of the form $(\alpha_p, \omega_{K,p} \beta_p)$ or
$(\omega_{K,p} \alpha_p, \beta_p)$ for unramified characters 
$\alpha_p$ and $\beta_p$ of $\Q_p^\times$.
\item If $p \in \Ram$ is inert in $L$, $\psi_p$ is either 
unramified or the product of $\theta_{p,L_p}$ or
$\theta_{p,L_p}^c$ and an unramified character.
\item If $p \in \Ram$ ramifies in $L$, $\psi_p$ is either 
unramified or the product of $\theta_{p,L_p}$ and an
unramified character.
\end{enumerate}
\end{lemma}

We can also explicate the equivalence relation 
$\Autf^1_K (n, \Delta; \pi') = \Autf^1_K (n, \Delta; \pi)$ 
for these representations.

\begin{lemma}
If
$\Autf^1_K (n, \Delta; \pi') = \Autf^1_K (n, \Delta; \pi)$
for $\pi$, $\pi' \in \Autf^1_{\Q} (n)$ which are 
automorphically induced from quadratic extensions,
they are necessarily induced from the same 
quadratic extension $L$. Furthermore,
for $\pi = {\rm AI}_{L/\Q} (\psi)$ and $\pi' = 
{\rm AI}_{L/\Q} (\psi')$ with $\psi$ and $\psi'$ as above,
the equivalence $\Autf^1_K (n, \Delta; \pi') = \Autf^1_K (n, \Delta; \pi)$
is true if and only if either $\delta = \psi' / 
\psi = \gamma \circ {\rm N}_{L / \Q}$
for some idele class character $\gamma$ of $\Q$ or $\delta = 
\psi' / \psi$ satisfies $\delta / \delta^c = \omega_K \circ {\rm N}_{L / \Q}$.
\end{lemma}

With these descriptions in hand, one can 
obtain a preliminary formula for the cardinality
of $\Autf^{1,bc}_K (n, \Delta)$,
which will in a second step be refined to 
a completely explicit expression. Using (\ref{H1Gammaa}) we can 
then compute the contribution to the cohomology of each 
individual group $\SL(2,\a)$.

To simplify the notation, we need the following definition: 
for an integer $n$ and an imaginary quadratic field $L$ 
define $\nu_{L,n} \in \{ 0, 1 \}$ as follows:
\begin{enumerate}
\item If $L$ is the field $\Q(\sqrt{-3})$, set 
\[ \nu_{L,n} = 
\left\{ \begin{array}{ll} 1,  & {\rm if}\ n \equiv 2 \, (3), \\
0, & \mbox{otherwise}. \end{array} \right.
\]
\item If $L = \Q(i)$, set
\[ \nu_{L,n} = 
\left\{ \begin{array}{ll} 1,  & {\rm if}\ n \equiv 1 \, (2), \\
0, & \mbox{otherwise}. \end{array} \right.
\]
\item If $L$ is not one of the two exceptional fields, 
we simply set $\nu_{L,n} = 1$ for all $n$.
\end{enumerate}

\begin{proposition} \label{cardbcrep}
\begin{enumerate} \item The cardinality of the set 
$\Autf^{1,bc}_K (n, \Delta)$ is given by 
\[
\frac{\left| \Autf^{1,bc}_K (n, \Delta) 
\right|}{\left| X(\Delta) \right|}
= \sum_{\pi \in \Autf^1_\Q (n)} 2^{-|R(\pi)|}
- \nu_{K,n} \frac{h_K}{2} - \sum_{L \in {\cal L} (K)} \nu_{L,n} 
2^{|\Ram| - |\Ram_L| - 2} h_L.
\]
\item For $L \in {\cal L} (K)$ we have  
\[
\frac{\left| \{ \Pi \in \Autf^{1,bc}_K (n, \Delta) \, | \,
\Pi \otimes 
\omega_L \circ {\rm N}_{K / \Q} \simeq \Pi
\} 
\right|}{\left| X(\Delta) \right|}
= \nu_{L,n} 2^{|\Ram| - |\Ram_L| - 2} h_L.
\]
\end{enumerate}
\end{proposition}

Combining this Proposition with (\ref{H1Gammaa}) and the fact
that $\left| \Autf^1_\Q (n; K) \right| = \nu_{K,n} 2^{|\Ram|-1} h_K$, we can immediately deduce:
\begin{proposition} \label{prelimformula}
The dimension of the base change part of the 
cohomology of the group $\SL(2,\a)$ is given by 
\begin{multline*}
\dim H_{\rm bc}^1 (\SL(2,\a), E_n) \\
=  \sum_{\pi \in \Autf^1_\Q (n)} 2^{-|R(\pi)|}
- \nu_{K,n} \frac{h_K}{2} - \sum_{L \in {\cal L} (K), \, 
\omega_L ({\rm N} (\a)) = -1} \nu_{L,n} 2^{|\Ram| - |\Ram_L| - 1} h_L.
\end{multline*}
\end{proposition}

\subsection{CM classes}\label{cmform}

As a consequence we also obtain the following results on cohomology spaces associated to CM automorphic forms. We introduce the following notation.

\begin{definition}
If $A_{\rm CM} \subseteq \Autf^1_K (n, \Delta)$ is the subset of
all automorphic representations automorphically induced from quadratic extensions of $K$, define
\[
H^1_{\rm CM} (\SL(2,\a), E_n) := H_{{\rm cusp},A_{\rm CM}}^1 (\SL(2,\a), E_n) \subseteq H_{{\rm cusp}}^1 (\SL(2,\a), E_n).
\]
\end{definition}
Note that the corresponding space is again independent of the choice of $\Delta$ with
$\Delta \cap \O^* = \{ 1 \}$.

First consider the intersection of this space with $H^1_{\rm bc}$. The following proposition follows immediately from (\ref{H1Gammaa}) and Proposition \ref{cardbcrep}, Part 2:

\begin{proposition} For $L \in {\cal L} (K)$ 
the 
representations in $\Autf^{1,bc}_K (n,\Delta)$ 
automorphically induced
from $KL$ 
contribute to $H^1_{\rm bc} (\SL(2,\a),E_n)$
a space of dimension
\[
\left\{ \begin{array}{ll} \nu_{L,n} 2^{|\Ram| - |\Ram_L| - 1} h_L,  &  
{\rm if}\ \omega_L ({\rm N} (\a)) = 1, \\
0, & \mbox{\rm otherwise}. \end{array} \right.
\]
\end{proposition}

We can also consider 
all
representations in 
$\Autf^1_K (n, \Delta)$ automorphically induced from
a fixed quadratic extension of $K$, necessarily of the form 
$KL$ for an imaginary quadratic extension $L$ as above. 
For this, let $L'$ be the real quadratic subfield of 
$LK$ and $h^+_{L'}$ its narrow ideal class number. The total number of such representations is then equal to
\begin{equation} \label{totalCM}
\frac{|X(\Delta)|}{2} \nu_{L,n} h_L h^{+}_{L'}.
\end{equation}
Consequently, we obtain:
\begin{proposition} For $L \in {\cal L} (K)$ 
the contribution of representations automorphically induced from $KL$
to $H_{\rm cusp}^1 (\SL(2,\a), E_n)$ 
has dimension
\[
\left\{ \begin{array}{ll} \nu_{L,n} h^{+}_{L'}  h_L,  &  {\rm if}\
\omega_L ({\rm N} (\a)) = 1, \\
0, & \mbox{\rm otherwise}. \end{array} \right.
\]
\end{proposition}

Note that this is precisely the contribution of twisted base 
change representations of the corresponding type times a factor of
$h^{+}_{L'} / 2^{|\Ram|-|\Ram_L|-1}$, which is the number 
of narrow ideal classes in a narrow genus of $L'$. 

The following relation between the dimension of the cohomology spaces for
 $\SL (2,\a)$ and $\SL(2,\O)$ follows immediately from (\ref{totalCM}) and (\ref{H1Gammaa}), this time applied with $A = \Autf^1_K (n, \Delta)$: 
\begin{proposition}
For any fractional ideal $\a$ of $K$ we have
\begin{equation}
\dim H_{\rm cusp}^1 (\SL(2,\a), E_n)
= \dim H_{\rm cusp}^1 ({\rm SL}(2 ,\o), E_n) 
\mbox{} -
\sum_{L \in {\cal L} (K), \, \omega_L ({\rm N} (\a)) = -1} 
\nu_{L,n} h^{+}_{L'} h_L.
\end{equation}
\end{proposition}

\begin{corollary} \label{CMbcCorollary}
For $L \in {\cal L} (K)$ 
there exist 
representations in $\Autf^1 (n, \Delta)$
automorphically induced from $KL$,
which are not twisted base changes from $\Q$, if and only 
if the narrow ideal class number $h^{+}_{L'}$
of the real quadratic subfield $L'$ of $KL$ is bigger than the 
corresponding number $g^{+}_{L'} = 2^{|\Ram (L')|-1}$ of genera. In this
case, the contribution of these representations 
to the dimension of $H_{\rm cusp}^1 (\SL(2,\a), E_n)$ is independent of 
$n$, if $L$ is not one of the two exceptional fields
$\Q(\sqrt{-3})$ and $\Q(i)$ and constant on residue classes modulo 
$3$ and $2$, respectively,
in the two exceptional cases.

The existence of such representations (i.~e. the failure of
the relation $H^1_{\rm CM} \subseteq H^1_{\rm bc}$ for the field $K$)
is equivalent to the existence of a real quadratic field $L'$ with 
$h^{+}_{L'} > g^{+}_{L'}$ and
$KL'/K$ unramified (equivalently, $d_{L'}$ divides the discriminant $d_K$ and $d_{L'}$ and
$d_K/ d_{L'}$ are coprime). 
\end{corollary}

In the following table we give the real quadratic fields
$L'=\Q(\sqrt{D})$ with the five smallest discriminants, for which 
the criterion of Corollary \ref{CMbcCorollary} is satisfied.
\begin{center}
\begin{tabular}{|c|c|c|}
\hline
$D = d_{L'}$ & $g^+$ & $h^+$ \\
\hline
$136=8\cdot 17$ & 2 & 4\\
$145=5\cdot 29$ & 2 & 4\\
$205=5\cdot 41$ & 2 & 4\\
$221=13\cdot 17$ & 2 & 8\\
$229$ & 1 & 3\\
\hline
\end{tabular}
\end{center}

\subsection{Dimension formulas}\label{dimform}

We now deduce from the preliminary formula of 
Proposition \ref{prelimformula} a completely explicit dimension formula
for $H^1_{\rm bc}$.
For $p \in \Ram$ let $\nu_p$ be the exact power of 
$p$ dividing the discriminant of $K$. We have $\nu_p = 1$ for
$p \neq 2$ and $\nu_2 = 2$ or $3$.

For any integer $n$ set
\begin{equation}
\varepsilon_n = \left\{ \begin{array}{ll} 
\frac{(-1)^{n/2}}{4}, & {\rm if}\ n \equiv 0 \, (2), \\
0, & \mbox{\rm otherwise}, \end{array} \right.
\qquad
{\rm and}\qquad
\mu_n = \left\{ \begin{array}{ll} 0, & {\rm if}\ n \equiv 1 \, (3), \\
- \frac{1}{3}, & {\rm if}\ n \equiv 2 \, (3), \\
\frac{1}{3}, & {\rm if}\ n \equiv 0 \, (3). \end{array} \right.
\end{equation}
\begin{theorem}\label{dimtheo} 
Let $K=\Q(\sqrt{d})$ be an imaginary quadratic number field 
with ring of integer $\O$ and 
$n$ be a non-negative integer. We have
\begin{multline*}
\dim H^1_{\rm bc} ({\rm SL}(2,\o), E_n) =
\left( \frac{1}{24} \prod_{p \in \Ram} (p^{\nu_p}+1) + c_2 (-1)^{n+1} \right) 
(n+1)\\ 
 - \nu_{K,n} \frac{h_K}{2} - 2^{|\Ram|-2}
+ c_4 \varepsilon_{n+2} + c_3 \mu_{n+2} + \delta_{n,0}
\end{multline*}
where $\delta_{n,0}$ stands for the Kronecker delta symbol.
The constant $c_2$
is given by
\[
c_2 = \left\{ \begin{array}{ll} 2^{|\Ram|-4}, & 
\mbox{\rm if $p \equiv 1 \, (4)$ for all $p \in \Ram$, $p \neq 2$,} \\
0, & \mbox{\rm otherwise.} \end{array} \right.
\]
The constants $c_4$ and $c_3$ are given by
\[
c_4 = \left\{ \begin{array}{ll} 2^{|\Ram|}, & 
\mbox{\rm if $p \equiv 1$ or $3 \, (8)$ for all $p \in \Ram$,} \\
2^{|\Ram|-1}, & \mbox{\rm if $2 \in \Ram$ and $p \equiv 1$ or 
$3 \, (8)$ for all $p \in \Ram$, $p \neq 2$,} \\
0, & \mbox{\rm otherwise,} \end{array} \right.
\]
and
\[
c_3 = \left\{ \begin{array}{ll} 2^{|\Ram|-1}, & 
\mbox{\rm if $p^{\nu_p} \equiv 1 \, (3)$ for all $p \in \Ram$,} \\
2^{|\Ram|-2}, & \mbox{\rm if $3 \in \Ram$ and $p^{\nu_p} 
\equiv 1 \, (3)$ for all $p \in \Ram$, $p \neq 3$,} \\ 
0, & \mbox{\rm otherwise.} \end{array} \right.
\]
Furthermore, for any fractional ideal $\a$ of $K$ we have
$$
\dim H_{\rm bc}^1 (\SL(2,\a), E_n) =  
\dim H^1_{\rm bc} ({\rm SL}(2,\o), 
E_n)
 - \sum_{L \neq K, \, \omega_L ({\rm N} (\a)) = -1} 
\nu_{L,k} 2^{|\Ram| - |\Ram_L| - 1} h_L.
$$
\end{theorem}

It is interesting to compare the resulting lower bound for the 
dimension of the cohomology group
$H^1_{\rm cusp} ({\rm SL}_2 (\o), \C)$
with the lower bound obtained by
Rohlfs \cite{Rohlfs}. In the work of Kr\"amer \cite{Kraemer} a lower 
bound for the dimension 
$H^1_{\rm cusp} ({\rm SL}_2 (\o), \C)$ agreeing with the bound 
$\dim H^1_{\rm bc} ({\rm SL}_2 (\o), \C)$ of the above theorem
is derived by a different method.
 
Given Proposition \ref{prelimformula}, the proof of 
Theorem \ref{dimtheo} rests on the computation of spaces of holomorphic 
elliptic modular forms with fixed local components. We
summarize the ingredients necessary to carry out this task in the remaining part of this
subsection, while omitting some elementary computations.
The possible local components are given in 
Definition \ref{loccomp}. The dimension computation is based on the 
following Proposition. 
\begin{proposition} \label{dimQ}
Let $N \ge 1$ and $k \ge 2$ be integers and 
$\sigma$ a representation of $G_N = {\rm SL}(2,\Z / N \Z)$ 
such that $\sigma(-I_2)$ is the scalar $(-1)^k$. 
Let $U_N \subseteq G_N$ be the subgroup of all upper 
triangular unipotent elements and $S_3$ and $S_4$ the images in $G_N$ 
of elements of ${\rm SL}(2,\Z)$ of order $3$ and $4$, respectively. Then
\begin{multline*}
\dim \Hom_{G_N} (\sigma, S_k (\Gamma (N))) = 
\frac{k-1}{12} \dim \sigma - \frac{1}{2} \dim \sigma^{U_N} \\
\mbox{} + \varepsilon_k \tr \sigma (S_4) + \rho_k 
\tr \sigma (S_3) + \delta_{k,2} \dim  
\sigma^{G_N}.
\end{multline*}
\end{proposition}

It is not difficult to prove this Proposition using the 
description of ${\rm SL}(2,\Z)$ as an amalgamated product 
of $\langle - S_3 \rangle$ and $\langle S_4 \rangle$ and the 
Eichler-Shimura isomorphism. 
Of course, it is also a consequence of the trace formula.
By taking for $\sigma$ a representation induced from the 
Borel subgroup, one recovers the classical dimension 
formulas for the group $\Gamma_0 (N)$ with nebentype (cf. \cite{dim}). 

It remains to make explicit the representations of 
${\rm SL}(2,\Z / N \Z)$ corresponding to the local conditions of Definition
\ref{loccomp} and to compute the terms appearing in Proposition \ref{dimQ}. 
Actually, we will consider irreducible representations 
$\sigma$ of ${\rm GL}(2,\Z / N \Z)$ occurring in the 
automorphic representations in question with multiplicity one 
and use their restrictions to ${\rm SL}(2,\Z / N \Z)$. 
The representations $\sigma$ can be written as tensor products of 
representations $\sigma_p$ of ${\rm GL}(2,\Z / p^{\nu_p} \Z)$ 
for $p \in R (\pi)$. 
For the principal series representations of Definition \ref{loccomp}
the necessary computation of 
dimensions and character values is standard and 
we refer to \cite{dim}. For the convenience of the reader we 
repeat the results here. The dimension of the corresponding 
representation $\sigma_p$ is $p^{\nu_p - 1} (p+1)$. 
The dimension of the space of $U_N$-invariants equals $2$. 
The character values are given by 
\[
\tr \sigma_p (S_3) = \left\{ \begin{array}{ll} 0, & 
{\rm if}\ p \equiv 2 \, (3), \\
1, & {\rm if}\ p = 3, \\
2, & {\rm if}\ p \equiv 1 \, (3), 
\end{array} \right.
\]
and
\[
\tr \sigma_p (S_4) = \left\{ \begin{array}{ll} 2 (-1)^{(p-1)/4}, & 
{\rm if}\ p \equiv 1 \, (4), \\
0, & \mbox{\rm otherwise}. \end{array} \right.
\]
The parity of $\sigma_p$ is equal to $\omega_{K,p} (-1)$.

For the supercuspidal components we can use the 
constructions of \cite{LocalLanglands}. 

\begin{lemma}
Let $p \in \Ram$ and $\Q_{p^2}$ be the unramified 
quadratic extension of $\Q_p$.
We can write the representation $\pi_p = {\rm AI}_{K_p / \Q_p} (\theta_p)$ 
as ${\rm AI}_{\Q_{p^2} / \Q_p} (\theta'_p)$ with a character
$\theta'_p$ of $\Q_{p^2}^\times$ fulfilling $\theta'_p / (\theta'_p)^\tau = 
\omega_{K,p} \circ N_{\Q_{p^2} / \Q_p}$. The minimal conductor of
such a character is $p^{\nu_p}$. Assume in the following 
that $\theta'_p$ has this minimal conductor. Then $\pi_p$
contains (with multiplicity one) a unique representation of 
${\rm GL}(2,\Z_p)$, and this representation factors through a 
representation $\sigma_p$ of ${\rm GL}(2,\Z / p^{\nu_p} \Z)$. 
The dimension of $\sigma_p$ is
$p^{\nu_p-1} (p-1)$. If $p$ is odd, $\sigma_p$ is the 
cuspidal representation of 
${\rm GL}(2,\mathbb{F}_p)$ associated to the character of 
$\mathbb{F}_{p^2}^\times$ obtained by restricting $\theta'_p$.

Furthermore, we have the following values for the 
traces at the torsion elements of ${\rm SL}(2,\Z)$:
\[
\tr \sigma_p (S_3) = \left\{ \begin{array}{ll} 0, & 
{\rm if}\ p \equiv 1 \, (3), \\
-2, & \mbox{{\rm if}\ \rm $p \equiv 2 \, (3)$, 
$p > 2$, or $p = 2$, $\nu_p = 3$,} \\
- 1, & {\rm if}\ p = 3, \\
2, & {\rm if}\ p = 2, \, \nu_p = 2,
\end{array} \right.
\]
and
\[
\tr \sigma_p (S_4) = \left\{ \begin{array}{ll} 2 (-1)^{(p-3)/4}, & 
{\rm if}\ p \equiv 3 \, (4), \\
0, & \mbox{\rm otherwise}. \end{array} \right.
\]
The parity of $\sigma_p$ is $- \omega_{K,p} (-1)$ except in the 
case $p = 2$, $\nu_p = 2$,
where it is $-1$ while $\omega_{K,p} (-1) = -1$.
\end{lemma} 

{\em Proof. } We quickly sketch the ingredients of the proof. Everything is based 
on the Tame Parametrization Theorem of \cite[Theorem 20.2]{LocalLanglands} 
and the explicit constructions in its proof. The case of odd $p$ is 
covered by [loc. cit., 19.1]. Dimensions and character values can then 
be read off from the standard description of cuspidal representations 
over finite fields in [loc. cit., 6.4]. For $p=2$ we need the 
constructions of [loc. cit., 19.3, 19.4] together with [loc. cit., 15.8] 
to describe the representations of ${\rm GL}(2,\Z_2)$ and to 
compute the character values. The dimension statement can 
be found in [loc. cit., Lemma 27.6]. 

We can now finish the proof of Theorem \ref{dimtheo}.
For any tensor product of local representations $\sigma_p$ the dimensions 
and character values are obtained by multiplication. The space of 
$U_N$-invariants is non-trivial only if all local components are principal 
series representations. The last term in Proposition \ref{dimQ} 
only appears for representations of level one (i.~e. for $\sigma$ 
the trivial representation). It remains to compute for each $n$ 
the sum of the contributions in Proposition \ref{prelimformula} for all 
possible combinations of local components with total parity $(-1)^n$. 
This is a tedious but elementary computation which we omit here.

\subsection{Bounds for the cohomology of Bianchi groups}\label{bound}

In this subsection we use the results of Section \ref{dimform} to give some 
bounds for the dimension of the cohomology spaces
$H^1(\SL(2,\O_{d}),E_n)$ as $|d|$ or $n$ 
go to infinity.
The first result is a more or less obvious consequence of Theorem 
\ref{dimtheo}.  
\begin{corollary}
Let $K$ be an imaginary quadratic number field 
with ring of integers $\O_K$. There is a bound $C_1 >0$ such that
$$\dim H^1(\SL(2,\O_K),E_n) \ge C_1 n \qquad ({\rm as}\ n \to \infty). $$
\end{corollary}
For the proof we only have to show that the coefficient of $n+1$ in the 
formula of Theorem \ref{dimtheo} is non-negative.   

The second result we see as a complement to the following theorem
which is proved in \cite{BGLS}.
\begin{theorem}
Let $G$ be a simple Lie group with Haar measure $\mu$. 
There is a constant $C_2 > 0$ 
such that $d(\G)$ is at most $C_2 {\rm vol}(G/\G)$  for every lattice $\G$ 
in $G$ where 
$d(\G)$ is the minimal number of generators of $\G$.
\end{theorem}
To use this theorem note that
$${\rm vol}\left(\SL(2,\C)/\SL(2,\O_{d})\right)
=\frac{|d|^{3/2}}{4\pi^2} \, \zeta_K(2)$$
where $\zeta_K(s)$ is the Dedekind zeta function of $K$, see 
\cite[Section 7]{EGM3}. It is easy to see that $\zeta_K(2)$ is bounded between 
two positive real numbers for all imaginary quadratic fields $K$. 
In view of Lemma \ref{triviesti} we obtain
\begin{corollary} Let $n$ be a (fixed) non-negative integer. 
There is a constant $C_3 > 0$ such that
$$\dim H^1(\SL(2,\O_{d}),E_n)\le C_3 |d|^{3/2}\qquad {as}\ |d|\to\infty.$$
\end{corollary}
We remark that this result also follows from the trace formula methods of Section \ref{dim} below.

Theorem \ref{dimtheo} implies
\begin{proposition} Let $n$ be a (fixed) non-negative integer.
 There is a constant $C_4 > 0$ such that
$$\dim H^1(\SL(2,\O_{d}),E_n)\ge C_4 |d|\qquad {as}\ |d|\to\infty.$$
\end{proposition}

\subsection{Base change and cocompact arithmetic groups}

By the work of Labesse-Schwermer \cite{LabesseSchwermer} and Rajan \cite{Rajan},
it is possible to use base change and the Jacquet-Langlands correspondence to study the cohomology of 
the cocompact arithmetic groups $\Gamma$ associated to quaternion algebras defined over fields $L$ such that the extension $L / L^{\rm tr}$, where $L^{\rm tr}$ is the maximal totally real subfield of $L$, is solvable (but not necessarily Galois). The resulting bound for $\dim H^1 (\Gamma, E_n)$ is determined by the dimension of certain spaces of Hilbert modular forms of weight $(n+2,2,\ldots,2)$ for $L^{\rm tr}$, and will be again linear on congruence classes. We do not go into the details here. Note that, in contrast to the case of the Bianchi groups, for a particular group $\Gamma$ the resulting bound will often be trivial. If we however consider the collection of all congruence subgroups, the conjecture of Waldhausen and Thurston has been verified for these groups by Rajan \cite{Rajan}. We can easily deduce from his arguments the following qualitative result:

\begin{proposition} Let $\Gamma$ be an arithmetic subgroup of ${\rm SL} (2, \mathbb{C})$ such that the field of definition $L$ of the corresponding quaternion algebra is a solvable extension of its maximal totally real subfield $L^{\rm tr}$. Then for every $c > 0$ there exists a finite index subgroup $\Delta$ of $\Gamma$ such that 
\[
\dim H^1 (\Delta, E_n) > c n
\]
for all $n \ge 0$.
\end{proposition}

{\em Proof.} By \cite{LabesseSchwermer,Rajan}, for a suitable $\Delta$ a lower bound for the dimension of the cohomology is given by the dimension of the space of Hilbert modular newforms of weight $(n+2,2,\ldots,2)$ for certain congruence subgroups of $\GL (2, L^{\rm tr})$. By adding additional local conditions, it is easily seen that it is possible to assume that the subgroups in question are torsion-free. Furthermore, their covolume can be made arbitrarily large by changing $\Delta$. Shimizu's dimension formula \cite{Shimizu} implies then that for any $c > 0$ we can find a subgroup $\Delta$ such that the dimension of the corresponding space of Hilbert modular forms is $\ge c (n+1)$ for all $n \ge 0$.  

If base change for ${\rm SL} (2)$ for arbitrary extensions of number fields was available, one could prove the corresponding result for all arithmetic lattices. This provides strong theoretical evidence for a positive answer to Question \ref{EQuestion} for arithmetic lattices.

\section{Upper bounds for the dimension of $H^1$}\label{dim}

In this section we derive upper bounds for the dimension of the cohomology spaces $H^1 (\Gamma, E_n)$ by using the generalized Eichler-Shimura isomorphism to transform the problem into a question on multiplicities of representations in $L^2 (\SL (2,\C) / \Gamma)$ and then using the trace formula to get information on these multiplicities. We first set up the form of the trace formula we need by specializing the work of W. Hoffmann \cite{HoffmannOrbitalIntegrals,HoffmannITF} to our situation. Then we consider the behavior of the dimension of $H^1 (\Gamma, E_n)$ as a function of $n$ and its behavior for fixed $n$ as $\Gamma$ varies over the standard congruence subgroups $\Gamma_0 (\mathfrak{a})$ of a Bianchi group (our result is in fact slightly more general, cf. Theorem \ref{DimCongruenceSubgroups} below).

\subsection{Review of the invariant trace formula for $\SL (2,\C)$}

Let $\Gamma$ be a general discrete subgroup of $G = {\rm SL}(2,\C)$ of finite 
covolume and consider the discrete part of $L^2 (G / \Gamma)$, which is a 
Hilbert space direct sum of irreducible unitary representations 
$\pi$ of $G$, each one occurring with a finite multiplicity $m(\pi, \Gamma)$. 
The irreducible unitary representations of $G$ most important to us are the principal series representations $\pi_{m,i\nu}$ for integers $m$ and real parameters $\nu$, which are obtained by unitary induction from the characters
\[
\sigma_{m,i \nu} (e^{u + i \theta}) = e^{i (\nu u + m \theta)} 
\]
of the maximal torus $T \simeq \C^\times$ of $G$. The representations $\pi_{m,i\nu}$ and
$\pi_{-m,-i\nu}$ are equivalent. 

We are interested in bounding the multiplicities $m (\pi_{m,0})$ from above. As explicated above, the dimension of $H^1_{\rm cusp} (\Gamma,E_n)$ is the same as the 
multiplicity $m (\pi_{2n+2,0})$.

We first recall the trace formula for $L^2 (G / \Gamma)$ in the form in which it has been explicitly worked out by Hoffmann for lattices of rank one \cite{HoffmannITF}. We specialize his results to the simpler case of $G = {\rm SL} (2, \C)$ and the trivial Hecke operator. As a preparation, we need to recall the basic relations between orbital integrals and principal series characters and the explicit form of the Plancherel formula, for which we use \cite[Ch. XI]{KnappBook} as a reference. We normalize measures as in [loc. cit.], i.e. we use the Haar measure on $G$ given by the product measure $dk \, dn \, da$ associated to the Iwasawa decomposition $G = K N A$, where $dk$ gives $K = {\rm SU} (2)$ total measure one, $dn$ is the standard measure on $N \simeq \C$, the upper triangular unipotent subgroup, and the measure $da$ on $A \simeq \mathbb{R}^{> 0}$ is $du$ in the parametrization $u \mapsto {\rm diag} \, (e^u, e^{-u})$.
We consider compactly supported functions $f \in C^\infty_c (G)$. To such a function is 
associated the function $F^T_f \in C^\infty_c (T)$ defined by
\[
F^T_f (t) = e^{2u} \int_{K \times N} f (k t n k^{-1}) dk dn.
\]
For $g \in G$ set 
\[
D_G (g) = \det_{\mathfrak{g} / \mathfrak{g}_g} (1 - {\rm ad} \, (g)).
\]
Then $|D_G(g)|^{1/2} = |t - t^{-1}|^2$ for a regular semisimple element $g$ with eigenvalues $t$ and $t^{-1}$, and $D_G (g) = 1$ otherwise.
For $g \in G$ define the orbital integral
\[
J_G (g, f) = |D_G (g)|^{1/2} \int_{G/G_g} f (x g x^{-1}) dx.
\]
Then $J_G (t,f) = F^T_f (t)$ for all regular elements $t \in T$ [loc. cit., (11.13), (11.14)]. From this one sees immediately that $F^T_f (t) = F^T_f (t^{-1})$. It is also easy to see that $F^T_f (\pm 1) = 8 \pi J_G (\pm n_1, f)$ for 
$n_1 = \left( \begin{array}{cc} 1 & 1 \\ 0 & 1 \end{array} \right)$.
The Fourier transform of $F^T_f$ yields the characters of the principal series representations:
\[
\Theta_{m,i\nu} (f) = \frac{1}{2 \pi} \int_T F^T_f (e^{u+i\theta}) e^{i (\nu u + m \theta)} du d\theta
\]
and
\[
F^T_f (e^{u + i \theta}) = \frac{1}{2 \pi} 
\sum_{m \in \Z} \int_{-\infty}^\infty
\Theta_{m,i\nu} (f) e^{-i(\nu u + m \theta)} d \nu.
\]
The Plancherel formula for $G$ is given by [loc. cit., Theorem 11.2] (up to a minor correction):
\begin{equation} \label{plancherel}
f(1) = \frac{1}{16 \pi^2} \sum_{m \in \Z} \int_{-\infty}^\infty
\Theta_{m,i\nu} (f) (m^2 + \nu^2) d \nu.
\end{equation}

For a discrete subgroup $\Gamma$ of $G$ of finite covolume let $\cal C$ 
be the set of all cuspidal parabolic subgroups of $G$, i.~e. 
of all parabolic subgroups fixing a cusp of $\Gamma$. 
Let $\Gamma (*)$ be the set of all semisimple elements of 
$\Gamma$ which do not fix a cusp together with the elements 
of $\Gamma \cap \{ \pm 1 \}$, and on the other hand $\Gamma_{ce}$ 
the set of semisimple elements of $\Gamma$ different from 
$\pm 1$ and stabilizing a cusp. For $\xi \in \Gamma_{ce}$ 
let $A(\xi)$ be the unique conjugate of the real 
torus $A \subseteq T$ in the centralizer of $\xi$. For the definition 
of the weight factor $v_\xi$ for $\xi \in \Gamma_{ce}$ 
we refer to \cite[p. 105]{HoffmannITF}. For each $P \in {\cal C}$ let $\Gamma_M (P)$ be the set of projections to a Levi component $L$ of $P$ of the elements of $\Gamma \cap P$. 
We define constants $C (P, \eta, \Gamma)$ 
(called $C_P (\eta n_1, \chi_\Gamma)$ 
in [loc. cit., p. 106]) in terms of 
Epstein zeta functions associated to $\eta \Gamma \cap N$. 
Namely, $C (P,\eta,\Gamma)$ is the constant term in the 
Laurent expansion at $z = 1$ of the function
\[
C (P, \eta, \Gamma; z) = \frac{2 {\rm vol} (N / 
\Gamma \cap N)}{|\Gamma_M (P)|} 
\sum_{\xi \in \eta \Gamma \cap N, \, \xi \neq 1} \frac{1}{| u(\xi)|^{2z}},
\]
where we choose $k_P \in K$ such that $k_P^{-1} P k_P$ is the standard upper triangular Borel subgroup $P_0$ and write
\[
\xi = k_P \left( \begin{array}{cc} 1 & u (\xi) \\ 0 & 1 \end{array} \right) k_P^{-1}, \quad \xi \in N.
\]
The absolute value of $u(\xi)$ does not depend on the choice of 
$k_P$. For $\eta = 1$ we can write $C(P,1,\Gamma) = 2 \pi 
\kappa_{\Lambda (P)} / |\Gamma_M (P)|$ for the lattice $\Lambda (P) = 
u (\Gamma \cap N)$ in $\C$, where $\kappa_\Lambda$ denotes the 
constant term in the expansion of $\frac{|\Lambda|}{\pi} 
\sum_{\lambda \in \Lambda \setminus \{ 0 \}} |\lambda|^{-2z}$ 
at $z = 1$ (this notation agrees with \cite[Lemma 6.5.2]{EGM3}).
We also need distributions $I_L (\eta)$, 
which are Arthur's invariant modifications of weighted orbital integrals (cf. [loc. cit., Sect. 5]).
Finally, let $\Phi (\sigma_{m,s})$ be the scattering matrix of $\Gamma$ defined in [loc. cit., p. 122] (and denoted by $S (\chi_\Gamma, \tilde{w}, \sigma_\Lambda)$ there) and $\phi = \det \Phi$ its determinant (with respect to a suitable identification of the vector spaces in question, the choice of which is unimportant). We can now quote \cite[Theorem 6.4]{HoffmannITF}, specialized to our situation.

\begin{theorem} For $f \in C^\infty_c (G)$ the trace of the corresponding convolution operator on the discrete part of $L^2 (G / \Gamma)$ is given by
\begin{eqnarray*}
\tr \pi_\Gamma^{\rm disc} (f) & = & 
\sum_{\{ \xi \}_{\Gamma} \subset \Gamma (*) }
{\rm vol} (G_\xi / \Gamma_\xi) | D_G (\xi) |^{-1/2} J_G (\xi, f) \\
& & \mbox{} +  \sum_{\{ \xi \}_{\Gamma} \subset \Gamma_{ce} }
{\rm vol} (G_\xi / \Gamma_\xi A (\xi)) v_\xi | D_G (\xi) |^{-1/2} J_G (\xi, f) \\
& & \mbox{} + \sum_{P \in {\cal C}, \, \eta \in \Gamma_M (P) \cap \{ \pm 1 \}} C (P, \eta, \Gamma) J_G (\eta \left( \begin{array}{cc} 1 & 1 \\ 0 & 1 \end{array} \right), f) \\
& & + \frac{1}{2} \sum_{P \in {\cal C}, \, \eta \in \Gamma_M (P)} | \Gamma_M (P) |^{-1} I_L (\eta, f) \\
& & + \frac{1}{4\pi} \sum_{m \in \Z} \int_{-\infty}^\infty
\frac{\phi' (\sigma_{m,i\nu})}{\phi (\sigma_{m,i\nu})} \Theta_{m,i\nu} (f) d \nu \\
& & - \frac{1}{4} \tr \Phi (\sigma_{0,0}) \Theta_{0,0} (f).
\end{eqnarray*}

\end{theorem}

The distributions $I_L (\eta)$ can be explicitly described 
in terms of the character values 
$\Theta_{m, i \nu} (f)$. We need here only the limiting case $\eta = 1$. 
Denote by $\psi (s) = \Gamma' (s) / \Gamma (s)$ the logarithmic derivative of the
gamma function. The following Proposition follows easily from Hoffmann's work in
\cite{HoffmannOrbitalIntegrals}.

\begin{proposition}
For the trivial element of $G$ the distribution $I_L (1)$ is given by  
\[
I_L (1, f) = \frac{1}{2 \pi} \sum_{m \in \Z} \int_{-\infty}^\infty 
\Omega_L (1, \sigma_{m,i \nu}) \Theta_{m,i \nu} (f) d \nu + 
\frac{1}{2} \Theta_{0,0} (f)
\]
with the function
\[
\Omega_L (1, \sigma_{m,i\nu}) = \psi (1) - {\rm Re} \, \psi (\frac{m+i\nu}{2}).
\]
\end{proposition}

Note that although the function $\psi$ has a simple pole at $s = 0$, 
the real part of $\psi (i \nu / 2)$ is continuous at $\nu = 0$, and in fact
${\rm Re} \, \psi (i \nu / 2) = {\rm Re} \, \psi (1 + i \nu / 2)$. 

{\em Proof.} Hoffmann considers invariant distributions $I_P$ closely related to $I_L$, cf. \cite[p. 58 bottom]{HoffmannOrbitalIntegrals} for their precise relation. The normalization factor $r_{\bar{P}P} (\sigma_{m,i\nu})$ there is up to a constant equal to $1 / (|m| + i \nu)$ (cf. \cite{KnappStein}). The distributions $I_P (1)$ are explicitly given by [loc. cit., p. 96, Corollary]. Note that we have only two roots $\alpha$ and $\bar{\alpha}$ and have to insert $\lambda (H_\alpha) = (m + i \nu) / 2$ and $\lambda (H_{\bar{\alpha}}) = (-m + i \nu) / 2$ into the expression given there.  Putting everything together and using the well-known relation $\Gamma (s) \Gamma (1-s) = \pi / \sin \pi s$, one obtains the formula above.

The reader may compare the resulting explicit trace formula, 
which involves only the function $F^T_f$ on $T$ and its 
Fourier transform, with the trace formula for 
$K$-biinvariant functions $f$ given in 
\cite[Theorem 6.5.1]{EGM3}. 

\subsection{The dimension of $H^1$: behavior with $n$}\label{dimn}

We now turn to the behavior of the multiplicities $m (\pi_{m,0}, \Gamma)$ as $m \to \infty$ for a fixed group $\Gamma$. The method extends to cohomological representations of real rank one groups which are not in the discrete series. It is an adaption of the method of \cite[Sect. 9]{DKVWeylLaw} for bounding the remainder term in Weyl's law. Our result is the following.

\begin{theorem} \label{DimVariationn}
For any discrete subgroup $\Gamma \subseteq G$ of finite covolume one has
\[
m (\pi_{m,0}) = O (m^2 / \log m), \quad m \to \infty.
\]
\end{theorem}

As an immediate consequence we have:

\begin{corollary}
For any discrete subgroup $\Gamma \subseteq G$ of finite covolume one has
\[
\dim H^1 (\Gamma, E_n) = O (n^2 / \log n), \quad n \to \infty.
\]
\end{corollary}

{\em Proof of the Theorem.} By passing to a finite index subgroup, we can assume that $\Gamma$ is torsion-free and $\Gamma_M (P) = \{ 1 \}$ for all $P$.

Let $m \ge 1$ and $g_0$ be an even $C^\infty$ function with support contained in $[-1,1]$, non-negative Fourier transform $h_0$ and $h_0 (0) > 0$. Consider the functions
\[
g (e^{u + i \theta}) = 2 \varepsilon g_0 (\varepsilon u) \cos m \theta
\]
on $T$, with $\varepsilon > 0$ being specified later. 
For $f \in C^\infty_c (G)$ with $F^T_f = g$ we have
\[
\Theta_{\pm m,i\nu} (f) = h_0 (\varepsilon^{-1} \nu),
\]
and $\Theta_{n,i\nu} (f) = 0$ for $|n| \neq m$. Insert $f$ into the trace formula and note
that because of our assumption on $\Gamma$, the sum in the second line is 
empty while the sums in the third and fourth line involve only $\eta = 1$. Also, the expression in the last line vanishes. Moving the integral involving the scattering matrix to the other side, we obtain an expression for
\begin{equation} \label{LHS}
2 \sum_\nu m (\pi_{m,i\nu}) h_0 (\varepsilon^{-1} \nu)
- \frac{1}{2 \pi} 
\int_{-\infty}^\infty
\frac{\phi' (\sigma_{m,i\nu})}{\phi (\sigma_{m,i\nu})} 
h_0 (\varepsilon^{-1} \nu) d \nu 
\end{equation}
as the sum of the remaining terms on the right hand side (i.~e. the sum of the first, third and fourth line).
As usual, we split the sum in the first line
as
\[
{\rm vol} (G / \Gamma) f (1) + \sum_{\{ \xi \}_{\Gamma} \subset \Gamma (*), \, \xi \neq 1}
{\rm vol} (G_\xi / \Gamma_\xi) | D_G (\xi) |^{-1/2} J_G (\xi, f).
\]
By the Plancherel formula, we can express the first term as
\[
\frac{{\rm vol} (G / \Gamma)}{8 \pi^2} \int_{-\infty}^\infty
h_0 (\varepsilon^{-1} \nu) (m^2 + \nu^2) d \nu
= C_1 \varepsilon m^2 + C_2 \varepsilon^3 
\]
with constants $C_1$ and $C_2$ depending only on $\Gamma$ and $h_0$. By \cite[pp. 90-91]{DKVWeylLaw}, since the absolute values of the orbital integrals $J_G (\xi,f)$ are
bounded independently of $m$, the second term can be estimated by $C_3 e^{C_4 / \varepsilon}$. The third line of the trace formula is a constant multiple of $F^T_f (1)$, and therefore $C_5 \varepsilon$. To estimate the weighted orbital integrals, we use the standard approximation $\psi (s) = \log s + O (1)$ for ${\rm Re} \, s \ge \delta > 0$ to get 
\[
|\Omega_L (1, \sigma_{m,i\nu})| \le \frac{1}{2} \log (m^2 + \nu^2) + C_6
\le \log m + \frac{\nu^2}{2 m^2} + C_6.
\]
From this one obtains 
\[
|I_L (1, f)| \le (C_7 + C_8 \log m) \varepsilon + C_9 \frac{\varepsilon^3}{m^2}.
\]
Taking $\varepsilon = c / \log m$ with a suitable constant $c$, one may conclude that with this choice the expression (\ref{LHS}) is $O (m^2 / \log m)$ as $m \to \infty$. Now, for each $m$ the determinant of the scattering matrix may be written as a Hadamard product
\[
\phi (\sigma_{m,s}) = \phi (\sigma_{m,0}) q_m^s \prod_{\eta \in P_m} \frac{s + \bar{\eta}}{s - \eta}
\]
with positive real constants $q_m$ bounded from above, where the product runs over the set $P_m$ of poles of $\phi_{m,s}$, which all have negative real part. Note that for $m \ge 1$ the Eisenstein series and the scattering determinant cannot have a pole on the real axis, since the corresponding induced representations do not have any unitarizable subquotients. Taking logarithmic derivatives, one sees that
$\log q_m - \phi' (\sigma_{m,i\nu}) / \phi (\sigma_{m,i\nu})$ is positive real for all $\nu$. Since $h_0$ was assumed to be non-negative, (\ref{LHS}) is therefore up to a term 
going to zero with $m$ an upper bound for $m (\pi_{m,0}) h_0 (0)$. The Theorem follows.

We remark that for congruence subgroups of the Bianchi groups ${\rm SL} (2, \O_K)$, $K$ imaginary quadratic, standard estimates for the logarithmic derivatives of Hecke $L$-functions imply that the contribution from the continuous spectrum in (\ref{LHS}) is $O (\varepsilon \log m)$ if $\varepsilon$ goes to zero for $m \to \infty$, and that it is therefore bounded with our choice of $\varepsilon$.

\subsection{The dimension of $H^1$: congruence subgroups of Bianchi groups}\label{dim2}

We now consider finite index subgroups of the Bianchi groups 
${\rm SL} (2, \O_K)$. For any non-zero ideal $\a$ of 
$\mathcal{O}_K$ we have the classical congruence subgroup
\[
\Gamma_0 (\a) = \left\{ \left( \begin{array}{cc} a & b \\ c & d 
\end{array} \right) \in {\rm SL}(2, \mathcal{O}_K) 
\, \ \vrule \ \, c \in \a \right\}.
\]
The index of $\Gamma_0 (\a)$ in ${\rm SL}(2,\mathcal{O}_K)$ is 
given by the multiplicative function
\[
\iota (\a) = {\rm N} (\a) \prod_{\mathfrak{p} \, | \, \a} 
\left( 1 + \frac{1}{{\rm N} (\mathfrak{p})} \right).
\]
We also need the following subgroups closely related to the 
principal congruence subgroups:
\[
\tilde{\Gamma}(\a) = \left\{  \left( \begin{array}{cc} a & b \\ c & d 
\end{array} \right) \in {\rm SL}(2,\mathcal{O}_K) 
\, \ \vrule \ \, a \equiv d \ {\rm mod }\ \a, \ b, c \in \a \right\}.
\]
The following theorem gives a bound for the multiplicity of a 
representation $\pi_{m,0}$ in $L^2 (G / \Gamma \cap \Gamma_0 (\a))$, 
$\Gamma$ a subgroup of finite index in ${\rm SL} (2, \mathcal{O}_K)$, 
which improves the trivial bound $O (\iota (\a))$ by a logarithm. 

\begin{theorem} \label{DimCongruenceSubgroups}
Let $\Gamma$ be a subgroup of finite index in 
${\rm SL}(2,\mathcal{O}_K)$, $K$ imaginary quadratic. 
Then for any fixed $m \ge 1$ we have
\[
m (\pi_{m,0}, \Gamma \cap \Gamma_0 (\a)) = 
O \left(\frac{\iota (\a)}{\log {\rm N} (\a)}\right), 
\quad {\rm N} (\a) \to \infty.
\]
\end{theorem}

This theorem can be regarded as a quantitative variant of the limit multiplicity results of de George-Wallach \cite{dW}, L\"{u}ck and Savin \cite{Sa} (which however concern towers of normal subgroups). Note also that for $\a = a \mathcal{O}_K$, $a$ a positive integer, 
we can get by base change arguments a lower bound of the form 
$C a = C {\rm N} (\a)^{1/2}$. If $\a$ and its conjugate are 
relatively prime, there is no non-trivial lower bound known.

\begin{corollary} 
Let $\Gamma$ be a subgroup of finite index in 
${\rm SL}(2,\mathcal{O}_K)$. 
Then for any fixed $n \ge 0$ we have
\[
\dim H^1 (\Gamma \cap \Gamma_0 (\a), E_n) = 
O \left(\frac{\iota (\a)}{\log {\rm N} (\a)}\right), 
\quad {\rm N} (\a) \to \infty.
\]
\end{corollary}

{\em Proof.} The corresponding assertion for the cuspidal part is an immediate consequence of Theorem \ref{DimCongruenceSubgroups}. To bound the dimension of the non-cuspidal part use Lemma \ref{LemmaCusps} below.

The proof of Theorem \ref{DimCongruenceSubgroups} is again based on the trace formula. 
By passing to a finite index subgroup, we 
can assume that $\Gamma$ is torsion-free and $\Gamma_M (P) = \{ 1 \}$ for all $P$.
Let $\Delta \subseteq \Gamma$ be a subgroup of finite index. We may then 
write (cf. \cite{Corwin}) for every $f \in C^\infty_c (G)$ the spectral side
\[
\tr \pi_\Delta^{\rm disc} (f)
- \frac{1}{4\pi} \sum_{m \in \Z} \int_{-\infty}^\infty
\frac{\phi'_\Delta (\sigma_{m,i\nu})}{\phi_\Delta (\sigma_{m,i\nu})} \Theta_{m,i\nu} (f) d \nu \\
+ \frac{1}{4} \tr \Phi_\Delta (\sigma_{0,0}) \Theta_{0,0} (f)
\]
as the sum
\begin{multline*}
[\Gamma : \Delta] {\rm vol} (G / \Gamma) f (1) + 
\sum_{\{ \xi \}_{\Gamma} \subset \Gamma (*), \, \xi \neq 1}
c_\Delta (\xi) {\rm vol} (G_\xi / \Gamma_\xi) | D_G (\xi) |^{-1/2} J_G (\xi, f) \\
\mbox{} + \sum_{P \in {\cal C}_\Delta} C (P, 1, \Delta) J_G (\left( \begin{array}{cc} 1 & 1 \\ 0 & 1 \end{array} \right), f) + \frac{1}{2} |{\cal C}_\Delta | I_L (1, f),
\end{multline*}
where we set
\[
c_\Delta (\xi) = \left| \{ \gamma \in \Delta \backslash \Gamma \, | \, \gamma \xi \gamma^{-1} \in \Delta \} \right|.
\]

We now need a sequence of elementary lemmas to deal with the parabolic and hyperbolic contributions.
The following well-known lemma is used to bound the parabolic contribution.

\begin{lemma} \label{LemmaCusps}
Let $\kappa$ be the multiplicative function defined by 
\[
\kappa (\mathfrak{p}^k) = \left\{ \begin{array}{ll} {\rm N} (\prim)^{k/2} + {\rm N} (\prim)^{k/2 - 1}, & k \equiv 0 \, (2), \\
				  2 {\rm N} (\prim)^{(k-1)/2}, & k \equiv 1 \, (2). 
\end{array} \right.
\]
Then we have
\[
|{\cal C}_{\Gamma \cap \Gamma_0 (\a)}| \le \kappa (\a) |{\cal C}_\Gamma| \le \frac{\iota (\a)}{\sqrt{{\rm N} (\a)}} |{\cal C}_\Gamma|,
\]
and the first inequality is an equality for $\Gamma = {\rm SL} (2, \O_K)$.
\end{lemma}

To deal with the hyperbolic contribution, we need to consider first the numbers $c_{\Gamma \cap \Gamma_0 (\a)} (\xi)$. The following lemma follows easily from the definitions.

\begin{lemma} \label{LemmaC}
Let $\xi \in \Gamma$ and $\mathfrak{b}$ be the largest divisor of $\a$ such that $\xi \in \tilde{\Gamma} (\b)$. Then 
\[
c_{\Gamma \cap \Gamma_0 (\a)} (\xi) \le c (\a, \b) \le 2^{\nu (\a)} {\rm N} (\mathfrak{b}),
\]
where $\nu (\a)$ denotes the number of prime divisors of $\a$ and $c$ is defined by
extending 
\[
c (\prim^k, \prim^r) =
\left\{ \begin{array}{ll} 2 {\rm N} (\prim)^r, & r < k, \\
			{\rm N} (\prim)^{k} + {\rm N} (\prim)^{k-1}, & r \ge k, 
\end{array} \right.
\]
multiplicatively.
\end{lemma}

For any semisimple element $\gamma \in G$ let its 
norm ${\rm N} (\gamma) \ge 1$ be the maximum value of 
$|t|^2$ for the two eigenvalues $t$ of $\gamma$. 
We need to estimate the number of $\Gamma$-conjugacy 
classes of bounded norm which are contained in $\tilde{\Gamma} (\b)$. 
Such an estimate can be deduced from the following well-known lemma.

\begin{lemma} \label{LemmaBoundCC}
There is a constant $B$ depending only on $\Gamma$ such that every semisimple conjugacy class $\{ \gamma \}_\Gamma$ in $\Gamma$ with ${\rm N} (\gamma) \le T$ contains a representative $\gamma = \left( \begin{array}{cc} a & b \\ c & d \end{array} \right)$ with 
$|a|^2, \, |b|^2, \, |c|^2, \, |d|^2 \le B T$.
\end{lemma}

The crude estimate of the next lemma is an easy consequence. The reader may verify that a better estimate would not change the final result (apart from the constant implicit in the $O$).
 
\begin{lemma} \label{LemmaNumberCC} For every $\delta > 0$ there is a constant $C$ depending on $\Gamma$ and $\delta$, such that for all non-zero ideals $\b$ of $\mathfrak{o}_K$ the number of $\Gamma$-conjugacy classes in $\Gamma (*)$ with norm $\le T$ which are contained in the normal subgroup $\Gamma \cap \tilde{\Gamma} (\b)$ of $\Gamma$ is bounded by $C T^{2+\delta} {\rm N} (\b)^{-2}$.
\end{lemma}

{\em Proof.} Apply Lemma \ref{LemmaBoundCC} to see that each such conjugacy class has a representative 
$\gamma = \left( \begin{array}{cc} a & b \\ c & d \end{array} \right)$
with $|a|^2, \, |b|^2, \, |c|^2, \, |d|^2 \le B T$. Furthermore, $bc \neq 0$ since the conjugacy class was assumed to lie in $\Gamma(*)$. The number of possible pairs $(b,c)$ corresponding to elements of $\tilde{\Gamma} (\b)$ is therefore bounded by $C' T^2 {\rm N} (\b)^{-2}$ with $C'$ depending only on $K$. For each such pair the number of possible entries $a$ and $d$ with $ad=1+bc \neq 0$ is clearly bounded by $O (T^\delta)$ with a constant depending only on $\delta$. This proves the assertion.

{\em Proof of Theorem \ref{DimCongruenceSubgroups}.}
We take a test function $f \in C^\infty_c (G)$ depending on $m \ge 1$ and $\varepsilon$ as above. We fix $m$ and assume at first only that $\varepsilon$ is bounded. Then the identity contribution to the trace formula for $\Delta = \Gamma \cap \Gamma_0 (\a)$ is bounded by $C_1 \iota (\a) \varepsilon$. The lattices $\Lambda (P) = u (\Delta \cap N)$ appearing in the definition of $C (P, 1, \Delta) = 2 \pi \kappa_{\Lambda (P)}$ are all invariant under a fixed order of the field $K$, and belong therefore to finitely many classes up to multiplication by elements of $K^\ast$. Using that $\kappa_\Lambda + \log |\Lambda|$ is invariant under such homotheties, this implies that the constants $\kappa_{\Lambda (P)}$ are bounded by $C + \log {\rm N} (\a)$ for a constant $C$.
By Lemma \ref{LemmaCusps} the parabolic contribution is therefore bounded by
$C_2 \iota (\a) (\log {\rm N} (\a)) {\rm N} (\a)^{-1/2} \varepsilon$.

As for the contribution of classes in $\Gamma(*)$, the estimates of Lemma \ref{LemmaC} and Lemma \ref{LemmaNumberCC} show that it is bounded by
\[
\left( 2^{\nu (a)} \sum_{\b \, | \, \a} \frac{1}{{\rm N} (\b)} \right) C_3 e^{C_4 / \varepsilon} \le C_5 (\mu) {\rm N} (\a)^\mu e^{C_4 / \varepsilon}
\]
for any $\mu > 0$. Taking 
$\varepsilon = c / \log {\rm N} (\a)$ with a suitable constant 
$c$, we see that the geometric side of the trace formula is 
indeed $O (\iota (\a) / \log {\rm N} (\a))$. A positivity argument 
as in the proof of Theorem \ref{DimVariationn} yields the result.

\section{Computational results for Bianchi groups}\label{compu}

This section contains computational results on the dimensions of the
cohomology groups $H^1(\G,E_n)$ where $\G=\SL(2,\a)$ is one of the 
Bianchi groups of Section \ref{present} and $n$ a non-negative integer. 
We also consider certain congruence subgroups of $\SL (2,\O_{-1})$. 
The method of computation is explained in Section 
\ref{Haeins}. 
Often the resulting systems of linear equations 
turned out to be far to big to do computations over the rational 
numbers. In these cases we were able to use the lower estimates 
of Section \ref{dimform} to deduce the dimension of the 
solution space over the complex numbers from the 
dimension over various finite fields.

\subsection{Dimensions of $H^1(\SL(2,\O),E_n)$}\label{compueucl}

Let us start with a little table. In Table \ref{ta1} we have listed 
the dimension of the cohomology spaces $H^1(\SL(2,\O_d),E_n)$ for 
$d=-1,\, -2,\, -3,\, -7,\, -11$ and
$0\le n \le 15$. 
To compare these values with the
dimensions of the spaces of lifted forms given in  
Proposition \ref{liftintro} or more generally in Theorem \ref{dimtheo}, it 
is important to know the codimension of the cuspidal cohomology.
Using \cite[Th. 8, Cor. 1]{SerreSL2} (and the further information contained there),
it can be easily checked that 
\begin{equation}
\dim H^1(\SL(2,\O_K),E_n)-\dim H_{\rm cusp}^1(\SL(2,\O_K),E_n)
= \nu_{K,n} h_K
\end{equation}
for all imaginary quadratic fields $K$, using the notation introduced in Section \ref{changegeneral}.
 
\begin{table}
\begin{center}
\begin{tabular}{|c|c|c|c|c|c|}
\hline
$n$ & $d=-1$ & $d=-2$ & $d=-3$ & $d=-7$ & $d=-11$ \\
\hline
$0$ & $0$ & $1$ & $0$ & $1$ & $1$  \\
$1$ & 1 & 1 & 0 & 1 & 1  \\
$2$ & 0 & 1 & 1 & 1 & 2  \\
$3$ & 1 & 2 & 0 & 1 & 2  \\
$4$ & 0 & 1 & 0 & 2 & 2  \\
$5$ & 2 & 3 & 1 & 2 & 3  \\
\hline
$6$ & 0 & 2 & 1 & 2 & 4  \\
$7$ & 3 & 4 & 1 & 3 & 4  \\
$8$ & 0 & 2 & 1 & 3 & 4  \\
$9$ & 3 & 5 & 1 & 3 & 5  \\
$10$ & 1 & 3 & 2 & 4 & {\bf 8}  \\
\hline
$11$ & 4 & 6 & 2 & 4 & 6  \\
$12$ & 0 & 3 & 1 & {\bf 6} & 6  \\
$13$ & 5 & 7 & 2 & 5 & 7  \\
$14$ & 1 & 4 & 3 & 5 & 8  \\
$15$ & 5 & 8 & 2 & 5 & 8  \\
\hline
\end{tabular}
\caption{\it Dimensions of $H^1(\SL(2,\O),E_n)$}\label{ta1}
\end{center}
\end{table}
We see that the cuspidal cohomology consists only of lifted forms
except in the two cases marked in boldface. These two cases will be analyzed 
more closely below. As a result of some heavy computer calculations 
we can report the following 
results.
\begin{proposition}
For $d=-1,\, -2,\, -3,\, -7,\, -11$ and $r_d$ as given in  
\begin{center}
\begin{tabular}{|c|c|c|c|c|c|}
\hline
$d$  & $-1$ & $-2$ & $-3$ & $-7$ & $-11$\\
\hline
${r}_{d}$ & $104$ & $141$ & $116$ & $132$ & $153$\\
\hline
\end{tabular}
\end{center}
we have 
$$H_{\rm cusp}^1(\SL(2,\O_d),E_n)=H_{\rm bc}^1(\SL(2,\O_d),E_n)$$
in the range $0\le n\le r_d$, except in the cases $d=-7$ and $n=12$, $d=-11$
and $n=10$, where $H^1_{\rm bc}$ has codimension two in $H^1_{\rm cusp}$.   
\end{proposition}
It remains to report the results of the computations in the 
non-euclidean cases. We have found the following:
\begin{proposition} Let $\G$ be one of the groups $\SL(2,\O_{d})$
with  $d=-19,\, -5,\, -6,\, -10,$ $ -14$ or 
$\SL(2,\a_{-5})$, $\SL(2,\a_{-6})$,
$\SL(2,\a_{-10})$, $\SL(2,\a_{-14})$, where the ideals 
$\a$ are as in Section \ref{present} and
let the non-negative integer $n$ be in the range $0\le n\le 60$. Then
$H_{\rm cusp}^1(\G,E_n)=H_{\rm bc}^1(\G,E_n)$. 
\end{proposition}

\subsection{Hecke operators on non-lifted cohomology classes}\label{hectab}

In this subsection we give the numerical values of some of the Hecke 
operators on the two spaces of non-lifted cohomology classes exhibited in 
Section \ref{compueucl}. 

\subsubsection{Hecke operators on $H^1(\SL(2,\O_{-7}),E_{12})$} 

We consider the prime element $\pi_{11}=2+\sqrt{-7}$
of $\O_{-7}$, which has degree one and norm $11$. 
By the methods described in Section \ref{hecke}
the characteristic polynomial of the corresponding
Hecke operator 
$$T_{\pi_{11}}: H^1(\SL(2,\O_{-7}),E_{12})\to  
H^1(\SL(2,\O_{-7}),E_{12})$$
can easily be computed to be
$$
P_{\pi_{11}}(X)=
\begin{array}{c}(X-9951764)(X^2+1877432\,X-54779120751344)\\ 
(X^3-2226532\,X^2-7410075237136\,X-1678794474022559168)
\end{array}.
$$
We know that there is a unique two-dimensional complement of the space of base change classes in the cohomology space $H^1_{\rm cusp} (\SL(2,\O_{-7}),E_{12})$ 
invariant under the Hecke operators. Identifying the 
Hecke operators on lifted classes (see Section \ref{hecke}) we infer 
that the kernel 
\begin{equation}
{\bf NL}(-7,12) = {\rm Ker}(T_{\pi_{11}}^2+1877432\,T_{\pi_{11}}-54779120751344)
\end{equation}
is this space of non-lifted classes. 
We write $L_\pi$ for the restriction of  the Hecke 
operators $T_\pi$, $\pi$ a prime element of $\O_{-7}$, to the space ${\bf NL}(-7,12)$. 
The following properties of the linear maps $L_\pi$ 
hold for all prime elements $\pi$ of $\O_{-7}$.
\begin{itemize}
\item[P7.1] $L_{-\pi}=-L_\pi$
\item[P7.2] $L_{\bar\pi}=-L_\pi^{\rm adj}$
\item[P7.3] If $\pi=p$ is a prime of degree two, then $L_{\pi}$ 
is an integer scalar denoted by $\lambda_p$.
\item[P7.4] After a suitable choice of basis for ${\bf NL}(-7,12)$, the
matrices giving the action of the Hecke operators $L_\pi$ have integral entries. 
\item[P7.5] The simultaneous splitting field for the Hecke 
operators $L_\pi$ is $\Q (\sqrt{7\cdot 239})$. 
\end{itemize}
Here $A^{\rm adj}$ stands for the adjoint of a linear map $A$. If $A$ 
is given by a two-by-two matrix, we have
$$\begin{pmatrix} a & b\\ c & d
\end{pmatrix}^{\rm adj}=
\begin{pmatrix} d & -b\\ -c & a
\end{pmatrix}.
$$
Property $P7.2$ follows by comparing the actions of the Hecke operator 
$T_\pi$ on the two cohomology spaces 
$$H^1(\SL(2,\O_{-7}),{\rm Sym}^{12}\otimes \overline{{\rm Sym}}^{12}),\qquad
H^1(\SL(2,\O_{-7}),\overline{{\rm Sym}}^{12}\otimes {\rm Sym}^{12}).$$
Property $P7.1$ is proved by computing the automorphism
$\epsilon$ induced by a matrix $E\in {\rm GL}(2,\O_{-7})$ of determinant $-1$
on $H^1(\PSL(2,\O_{-7}),E_{12})$. 
This property implies that no non-zero class in ${\bf NL}(-7,12)$ 
is the restriction of a cohomology class in $H^1({\rm GL}(2,\O_{-7}),E_{12})$. 
The rest of the above properties is clear.

Examples of the scalars $\lambda_p$ for primes of degree two are contained in
Table \ref{tahe71}. Examples of the integral matrices corresponding to the 
linear maps $L_\pi$ for primes $\pi$ of degree one  
are contained in Tables \ref{tahe72}, \ref{tahe73} and \ref{tahe74}. 

\begin{table}
\begin{center}
\begin{tabular}{|c|c|}
\hline
$p$  & $\lambda_p$ \\
\hline
$3$ & $-1939626$\\
$5$ & -747491750\\
$13$ & -252803502896086\\
$17$ & 4756247617499746 \\
$19$ & 5094169624293878\\
$31$ & -30279773153264109058\\
$41$ & -948454707467278569518\\
$47$ & -9168990821180522751074\\
$59$ & 123833654051598471764998\\
$61$ & -105716258627702854298998\\
$73$ & -707186203752039245531566\\
$83$ & -5005894274852029376014346\\
$89$ &  -980936263375178621227022\\
$97$ & 84206314563458516168628866\\
\hline
\end{tabular}
\caption{\it Scalars $\lambda_p$ for $L_p$ on ${\bf NL}(-7,12)$}\label{tahe71}
\end{center}
\end{table}

\begin{table}
\begin{center}
\begin{tabular}{|c|c|c||c|c|c|}
\hline
$p$  & $\pi$ & $A_\pi$ & $p$ & $\pi$ & $A_\pi$ 
\\
\hline
$2$ & $\omega$ & 
$\begin{pmatrix} 0 & 1 \\ 14432 & -50 \end{pmatrix}$ &
$7$ & $-1+2\omega$ & 
$\begin{pmatrix} 44800  &   1792\\
25862144 &  -44800
 \end{pmatrix}$ \\
\hline
$11$ & $1+2\omega$ & 
$\begin{pmatrix} 581284  &   60800\\
877465600 & -2458716\end{pmatrix}$ &
$23$ & $3+2\omega$ & 
$\begin{pmatrix} -257854600 &    4457728\\
64333930496 & -480741000 \end{pmatrix}$\\ 
\hline
\end{tabular}
\caption{\it Hecke operators $L_\pi$ on ${\bf NL}(-7,12)$}\label{tahe72}
\end{center}
\end{table}

\begin{table}
\begin{center}
\begin{tabular}{|c|c|c|}
\hline
$p$  & $\pi$ & $A_\pi$\\
\hline
$29$ & $-1+4\omega$ & 
$\begin{pmatrix}  -114226222  &   -627200\\
-9051750400 &  -82866222 \end{pmatrix}$\\
\hline
 $37$ & $1+4\omega$ & 
$\begin{pmatrix}  8869653750 &    -61610496\\
-889162678272 &  11950178550 \end{pmatrix}$\\
\hline
$43$ & $5+2\omega$ & 
$\begin{pmatrix}
42167274700  &   293147008\\
4230697619456 &  27509924300
\end{pmatrix}$\\
\hline
$53$ & $3+4\omega$ &
$\begin{pmatrix}
 -229421381350  &    -843922944\\
-12179495927808 &  -187225234150
\end{pmatrix}$\\
\hline
$67$ & $-1+6\omega$ &
$\begin{pmatrix}
914163852100 &    -1805341824\\
-26054693203968 & 1004430943300
\end{pmatrix}$\\
\hline
$71$ & $7+2\omega$ &
$\begin{pmatrix}
-600257601424 &    -5497094400\\
-79334066380800 &  -325402881424
\end{pmatrix}$\\
\hline
\end{tabular}
\caption{\it Hecke operators $L_\pi$ on ${\bf NL}(-7,12)$}\label{tahe73}
\end{center}
\end{table}

\begin{table}
\begin{center}
\begin{tabular}{|c|c|c|}
\hline
$p$  & $\pi$ & $A_\pi$\\
\hline
$79$ & $1+6\omega$ &
$\begin{pmatrix}
-775382036248 &    -11236492800\\
-162165064089600 &   -213557396248
\end{pmatrix}$\\
\hline
$107$ & $9+2\omega$ &
$\begin{pmatrix}
11411424109300  &   66437695872\\
958828826824704 &  8089539315700
\end{pmatrix}$\\
\hline
$109$ & $7+4\omega$ &
$\begin{pmatrix}
8058528373122   &   78253401600\\
1129353091891200 &   4145858293122
\end{pmatrix}$\\
\hline
$113$ & $3-8\omega$ &
$\begin{pmatrix}
4624127056750 &    42202431488\\
609065491234816 &  2514005482350
\end{pmatrix}$\\
\hline
$127$ & $5+6\omega$ &
$\begin{pmatrix}
4687450108000  &   367392485376\\
5302208348946432 & -13682174160800
\end{pmatrix}$\\
\hline
$137$ & $1+8\omega$ &
$\begin{pmatrix}
   79180120345450  &   -113527447552\\
-1638428123070464  &  84856492723050
\end{pmatrix}$\\
\hline
$149$ & $9+4\omega$ &
$\begin{pmatrix}
-71276735378522  &   562548416000\\
8118698739712000 & -99404156178522
\end{pmatrix}$\\
\hline
\end{tabular}
\caption{\it Hecke operators $L_\pi$ on ${\bf NL}(-7,12)$}\label{tahe74}
\end{center}
\end{table}

\subsubsection{Hecke operators on $H^1(\SL(2,\O_{-11}),E_{10})$} 

We consider the prime element $\pi_{3} =(1+\sqrt{-11})/2$
of $\O_{-11}$, which has degree one and norm $3$. 
By the methods described in Section \ref{hecke}
the characteristic polynomial of the corresponding
Hecke operator 
$$T_{\pi_{3}}: H^1(\SL(2,\O_{-11}),E_{10})\to  
H^1(\SL(2,\O_{-11}),E_{10})$$
can easily be computed to be
$$
P_{\pi_{3}}(X)=
\begin{array}{c}
(X-252)(X-67)(X^2+700\,X+40671) \\
(X^4+403\,X^3-439713\,X^2-113276475\,X+1097145000)
\end{array}.
$$
We know that there is a unique two-dimensional complement of the space of base change classes in the cohomology space $H^1_{\rm cusp} (\SL(2,\O_{-11}),E_{10})$ 
invariant under the Hecke operators. We infer that the kernel 
\begin{equation}
{\bf NL}(-11,10)={\rm Ker}(T_{\pi_{3}}^2+700\,T_{\pi_{3}}+40671)
\end{equation}
is equal to this space of non-lifted classes. 
We write $L_\pi$ for the restriction of  the Hecke 
operator $T_\pi$, $\pi$ a prime element of $\O_{-11}$, 
to the space ${\bf NL}(-11,10)$. 
The following properties of the linear maps $L_\pi$ 
hold for all prime elements $\pi$ of $\O_{-11}$.
\begin{itemize}
\item[P11.1] $L_{-\pi}=L_\pi$
\item[P11.2] $L_{\bar\pi}=-L_\pi^{\rm adj}$
\item[P11.3] If $\pi=p$ is a prime of degree two, then $L_{\pi}$ 
is an integer scalar denoted by $\mu_p$.
\item[P11.4] After a suitable choice of basis for ${\bf NL}(-11,10)$, the
matrices giving the action of the Hecke operators $L_\pi$ have integral entries. 
\item[P11.5] The simultaneous splitting field for the Hecke 
operators $L_\pi$ is $\Q (\sqrt{11\cdot 43\cdot 173})$. 
\item[P11.6] 
The space ${\bf NL}(-11,10)$ 
is the restriction of a subspace of 
$H^1({\rm GL}(2,\O_{-11}),E_{10})$. 
\end{itemize}

The case $d=-11$, $k=10$, differs from the case $d=-7$, $k=12$, 
since we have $L_{-\pi}=L_\pi$ for $d=-11$ and
and $L_{-\pi}=-L_\pi$ for $d=-7$. In the case $d=-11$ this leads directly to 
property P11.6.

Examples of the scalars $\mu_p$ for primes of degree two are contained in
Table \ref{tahe111}. Examples of the integral matrices corresponding to the 
linear maps $L_\pi$ for primes $\pi$ of degree one  
are contained in Tables \ref{tahe112}, \ref{tahe113}, \ref{tahe114}. 

\begin{table}
\begin{center}
\begin{tabular}{|c|c|}
\hline
$p$  & $\mu_p$ \\
\hline
$2$ & $-80$\\
$7$ & $-818885550$\\
$13$ & $1235127129530$\\
$17$ & $45387811032610$\\
$19$ & $-95158947964038$\\
$29$ & $-8701360899198758$\\
$41$ & $-429545462511285518$\\
$43$ & $638559982027780650$\\
$61$ & $10654154103002912922$\\
$73$ & $34915634850910529970$\\
$79$ & $-688424011186184859358$\\
$83$ & $33668143605728046010$\\
$101$ & $15523742571431406528202$\\
\hline
\end{tabular}
\caption{\it Scalars $\mu_p$ for $L_p$ on ${\bf NL}(-11,10)$}\label{tahe111}
\end{center}
\end{table}

\begin{table}
\begin{center}
\begin{tabular}{|c|c|c||c|c|c|}
\hline
$p$  & $\pi$ & $A_\pi$ & $p$ & $\pi$ & $A_\pi$ 
\\
\hline
$3$ & $\omega$ & 
$\begin{pmatrix} 0 & 1 \\ -40671 & -700 \end{pmatrix}$ &
$5$ & $-2+\omega$ & 
$\begin{pmatrix} -14203 & -26 \\ 1057446 & 3997 \end{pmatrix}$ \\
\hline

$11$ & $-1+ 2\omega$ & 
$\begin{pmatrix} -117612 & 0 \\ 0 & -117612 \end{pmatrix}$ &
$23$ & $-5+\omega$ & 
$\begin{pmatrix} 44565050 & 22561 \\ -917578431 & 28772350 \end{pmatrix}$\\ 
\hline

$31$ & $-4+3\omega$ & 
$\begin{pmatrix} -124944582 & -577125 \\ 23472250875 & 279042918 \end{pmatrix}$
&
 $37$ & $-5+3\omega$ & 
$\begin{pmatrix} 351981325 & 819882 \\ -33345420822 & -221936075 \end{pmatrix}$
 \\
\hline
\end{tabular}
\caption{\it Hecke operators $L_\pi$ on ${\bf NL}(-11,10)$}\label{tahe112}
\end{center}
\end{table}

\begin{table}
\begin{center}
\begin{tabular}{|c|c|c|}
\hline
$p$  & $\pi$ & $A_\pi$\\
\hline

$47$ & $-7+2\omega$ & 
$\begin{pmatrix} 2959574800  &  5646848 \\ -229662955008  & -993218800 \end{pmatrix}$
 \\
\hline
$53$ & $-5+ 4\omega$ & 
$\begin{pmatrix} -3591316050 &  -868224 \\ 35311538304 & -2983559250 \end{pmatrix}$
 \\
\hline
$59$ & $-8+ \omega$ & 
$\begin{pmatrix} 2044859460 & 10611525  \\ -431581333275  & -5383208040 \end{pmatrix}$
 \\
\hline
$67$ & $-8+ 3\omega$ & 
$\begin{pmatrix} 5506303200 & -13041567  \\ 530413571457  & 14635400100 \end{pmatrix}$
 \\
\hline
$71$ & $-4+ 5\omega$ & 
$\begin{pmatrix}-20524885978  &  -34309625  \\  1395406758375 & 3491851522
\end{pmatrix}$
 \\
\hline
$89$ & $-7+ 5\omega$ & 
$\begin{pmatrix}   19167086435  &     60342700\\
-2454197951700 &  -23072803565 
\end{pmatrix}$
 \\

\hline
\end{tabular}
\caption{\it Hecke operators $L_\pi$ on ${\bf NL}(-11,10)$}\label{tahe113}
\end{center}
\end{table}

\begin{table}
\begin{center}
\begin{tabular}{|c|c|c|}
\hline
$p$  & $\pi$ & $A_\pi$\\
\hline
$97$ & $-10+ 3\omega$ & 
$\begin{pmatrix}  33903167375  &     94396752\\
-3839210300592 &  -32174559025 
\end{pmatrix}$
 \\
\hline
$103$ & $-5+ 6\omega$ & 
$\begin{pmatrix}  -127510128200 &     -435253824\\
17702208275904 &  177167548600 
\end{pmatrix}$
 \\
\hline
$113$ & $-11+ \omega$ & 
$\begin{pmatrix} -257969686425 &    -640536456\\
26051258201976 &  190405832775 
\end{pmatrix}$
 \\
\hline
$137$ & $-5+ 7\omega$ & 
$\begin{pmatrix} 500270562475  &   -968684668\\
39397374132228 & 1178349830075 
\end{pmatrix}$
 \\

\hline
$157$ & $-13+ 3\omega$ & 
$\begin{pmatrix}  -2300340926975 &     -6625408122\\
269461973729862 &  2337444758425
\end{pmatrix}$
 \\

\hline
$163$ & $-11+ 6\omega$ & 
$\begin{pmatrix} -174488742500  &     315270144\\
-12822352026624 &  -395177843300
\end{pmatrix}$
 \\
\hline
$179$ & $-13+ 5\omega$ & 
$\begin{pmatrix} -933096107380  &   -5594462825\\
227532397555575 &  2983027870120
\end{pmatrix}$
 \\

\hline
\end{tabular}
\caption{\it Hecke operators $L_\pi$ on ${\bf NL}(-11,10)$}\label{tahe114}
\end{center}
\end{table}

\subsection{Cohomology of congruence subgroups}\label{CongruenceComp}

In this subsection we give some computational results concerning the dimensions 
of the cohomology groups $H^1(\Gamma,E_n)$ where $\Gamma \subseteq \SL(2,\O_{-1})$ 
is a congruence subgroup.

\subsubsection{The case of trivial coefficients}\label{Ga0}

Here we consider 
the congruence subgroups 
$$\Gamma^0(\mathfrak{p})=\left\{ \begin{pmatrix} a & b\\ c & d \end{pmatrix} \in
\SL(2,\O_{-1})\ \vrule \ b\in \mathfrak{p}\, \right\}$$ of $\SL (2, \O_{-1})$,
where $\mathfrak{p}$ is a prime ideal of $\O = \O_{-1}$ of degree one. Note that 
$\Gamma^0 (\mathfrak{p})$ is conjugate in $\SL (2, \O)$ to the congruence subgroup $\Gamma_0 (\mathfrak{p})$ considered in Section \ref{dim2}.
The norm of $\mathfrak{p}$ is either equal to $2$ or a rational prime $p$ congruent 
to $1$ modulo $4$. The index of $\Gamma^0(\mathfrak{p})$ in $\SL(2,\O_{-1})$ is 
$p+1$. The cohomology groups $H^1(\Gamma^0(\mathfrak{p}),\C)$ are particularly 
interesting for number theory since their non-vanishing is conjectured to be related to the 
existence of certain elliptic curves (or more generally abelian varieties) defined over $K=\Q(i)$ (cf. \cite{C,GHM,GM}).
We shall report here on extensive computations of the dimensions of
the spaces $H^1(\Gamma^0(\mathfrak{p}),\C)$.
Note that we have $H^1 (\Gamma^0 (\mathfrak{p}), \C) = 
H^1_{\rm cusp} (\Gamma^0 (\mathfrak{p}), \C)$ and $H^1 (\SL (2,\O),\C) = 0$, and that
$H^1 (\Gamma^0 (\mathfrak{p}), \C)$ consists therefore entirely of new classes (cf. Section \ref{Gano}).

The elements $A^i$, $0\le i\le p-1$, and $B$ (cf. Section \ref{present})
form a system of coset representatives for $\Gamma^0(\mathfrak{p})$ in 
$\SL(2,\O_{-1})$.  
From this we obtain the following generating system for 
$\Gamma^0(\mathfrak{p})$:
$$A^p,\ BAB,\ BUB,\ UA^\rho, \ A^{i'}BA^i, \qquad 1\le i \le p-1,\ ii' \equiv 1 
\ {\rm mod}\ p,
$$ 
where $\rho^2+1=0$ in the field $\F_p$.
From the presentation (\ref{bia1}) we may compute a 
presentation of the finitely generated abelian group $\Gamma^0(\mathfrak{p})^{\rm ab}$ 
and in particular the dimension of $\Gamma^0(\mathfrak{p})^{\rm ab} \otimes \C$ 
(which is the same as the dimension of $H^1(\Gamma^0(\mathfrak{p}),\C)$). 
This computation may be speeded up in the following way, i.~e.
the presentation of $\Gamma^0(\mathfrak{p})$ obtained from the Reidemeister-Schreier method
can be simplified a lot. We shall describe a result contained in \cite{GHM} which gives such a 
simplification.    
  
Let $R$ be a commutative ring and let ${\bf P}^1(R,p)$ be a $(p+1)$-dimensional 
free $R$-module with basis  
$u_x$ indexed by the projective line $\PP^1(\F_p)$.
This module has rank $p+1$ and is a
${\rm PGL}(2,\F_p)$-module by the natural permutation action on the basis elements.
Let ${\bf U}(R,p)$ be the submodule of ${\bf P}^1(R,p)$ generated by $u_0$ 
and the elements 
\begin{equation}
u_x+u_{Bx},\ \ u_x+u_{Wx},\ \ u_x+u_{Sx}+u_{S^2x},\ \ u_x+u_{Yx}+u_{Y^2x},
\qquad x\in\PP^1(\F_p),   
\end{equation}
with the matrices
$$B =\begin{pmatrix} 0 & 1\\ -1 & 0\end{pmatrix},
\qquad W=\begin{pmatrix} 0 & 1\\ 1 & 0\end{pmatrix},
\qquad S=\begin{pmatrix} 1 & -1\\ 1 & 0\end{pmatrix},
\qquad Y=\begin{pmatrix} -\rho & 1\\ 1 & 0\end{pmatrix}.
$$
Define 
$\Phi_\mathfrak{p}: {\bf P}^1(R,p)\to \Gamma^0(\mathfrak{p})^{\rm ab} \otimes R$ by setting 
$\Phi_\mathfrak{p}(u_i) =A^{i'}BA^i$ for $i\in \F_p$ with $i\ne 0$ and  
$\Phi_\mathfrak{p}(u_0) = \Phi_\mathfrak{p}(u_\infty)=0$.  
The results of \cite[Section 3]{GHM} imply
that the map $\Phi_\mathfrak{p}$ is a surjective group homomorphism with kernel equal  
to ${\bf U}(R,p)$, if $R$ is a field of 
characteristic $0$. 

To avoid heavy integer computations, we take $R=\F_q$ for various 
(small) primes $q$, compute the dimension of 
${\bf P}^1(R,p)/{\bf U}(R,p)$ and set
\begin{equation}
\dim_{\le x}\, H^1(\Gamma^0(\mathfrak{p}),\C) 
=\inf\  \{\, \dim_{\F_q}{\bf P}^1(\F_q,p)/{\bf U}(\F_q,p) \,\}  
\end{equation}
where $q$ ranges over all primes below $x$. Of course, if this number is 
zero then also $H^1(\Gamma^0(\mathfrak{p}),\C) = 0$, and if
$x$ is sufficiently large $\dim_{\le x} H^1(\Gamma^0(\mathfrak{p}),\C)$ 
will be equal to the dimension of   
$H^1(\Gamma^0(\mathfrak{p}),\C)$.

In Table \ref{ta22} we give the norms 
of the degree one primes $\mathfrak{p}$ in $\O_{-1}$ with $N(\mathfrak{p})\le 20000$ and  
$\dim_{\le 500}\, H^1(\Gamma^0(\mathfrak{p}),\C)=1$.
Table \ref{ta23} covers the same range and gives the 
norms of the degree one primes $\mathfrak{p}$  
with $\dim_{\le 500}\, H^1(\Gamma^0(\mathfrak{p}),\C)=2$.
The norms of the primes with $\dim_{\le 500}\, H^1(\Gamma^0(\mathfrak{p}),\C)=3$ 
are
$941$, $1777$, $5113$. 
Those with $\dim_{\le 500}=4$ are 
$8893$, $17021$. The values $5$ and $6$ are attained for 
$4517$, $5309$ respectively. There is no prime $\mathfrak{p}$ with 
$N(\mathfrak{p})\le 20000$ and $\dim_{\le 500} \ge 7$. 

\begin{table}
\begin{center}
\begin{tabular}{|cccccccc|}
\hline

137 &
233 &
257 &
277 &
509 &
569 &
733 & 
977 
\\

1009  &   
1013  &
1021 &
1049   & 
1153 &
1277 &
1373 &
1489

\\

1493 &
1753 & 
1997 & 
2053 & 
2081 & 
2377 & 
2441 &
2521

\\

2609 & 
2729 & 
2917 & 
3109 &
3361 &
3929 & 
4013 & 
4177

\\

4289 &
4421 &
4597 &
4621 &
4721 &
5021 &
5237 & 
5741 

\\

5749 &
5801 &
6029 &
6361 &
6701 &
6781 &
6793 &
6857
\\

6949 &
7001 &
7069 &
7121 &
7793 &
7937 &
8297 &
8377
\\

8461 &
8513 &
8537 &
8753 &
9041 &
9413 &
10357 &
10369 
\\
10477 &
10657 &
10729 &
10861 &
10937 &
11701 &
11953 &
12253
\\
12553 &
13381 &
13457 &
13633 &
15161 &
15497 &
15569 &
15629
\\
15749 &
16097 &
16349 &
16649 &
16673 &
17209 &
17921 &
18289 
\\
18553 &
18701 &
18869 &
18913 &
19213 &
19417 &
19841 &
19997 \\
\hline
\end{tabular}
\caption{\it Norms of degree one primes $\mathfrak{p}$ in $\O_{-1}$ 
with $\dim_{\le 500}\, H^1(\Gamma^0(\mathfrak{p}),\C)=1$}\label{ta22}
\end{center}
\end{table}

\begin{table}
\begin{center}
\begin{tabular}{|cccccccc|}
\hline

433 &
709 &
757 & 
853 &
953 &
1321 &
1549 &
1901 
\\

1973 &
2657 &
2753 &
3313 &
3469 & 
3529 &
3637 &
3877 
\\

5849 &
5857 &
6689 &
7577 &
8081 &
9349 &
9629 &
11437
\\

12269 &
12953 &
13093 &
13477 &
15761 &
16921 & 
17033 &
18757 
\\
19237 &
19937 & & & & & & 
\\
\hline
\end{tabular}
\caption{\it Norms of degree one primes $\mathfrak{p}$ in $\O_{-1}$ 
with $\dim_{\le 500}\, H^1(\Gamma^0(\mathfrak{p}),\C)=2$}\label{ta23}
\end{center}
\end{table}

In an even more extensive search we have gone 
through the degree one primes $\mathfrak{p}$ in $\O_{-1}$ with $N(\mathfrak{p})\le 60000$ and  
have computed $\dim_{\le 500}\, H^1(\Gamma^0(\mathfrak{p}),\C)$.
There are altogether $3018$ primes below $60000$ which are congruent 
to $1$ modulo $4$.
In the following table we give the number $N(r,60000)$ of 
such primes with $\dim_{\le 500}\, H^1(\Gamma^0(\mathfrak{p}),\C)=r$.

\begin{center}
\begin{tabular}{|c|c|c|c|c|c|c|c|c|c|c|}
\hline
r  & 0 & 1 & 2 & 3 & 4 & 5 & 6 & 7 & 8 & $\ge 9$\\
\hline
$N(r,60000)$ & 2728 & 198 & 73 & 11 & 4 & 1 & 1 & 1 & 1 & 0\\
\hline
\end{tabular}
\end{center}

\smallskip
The value $8$ is attained for the prime $58313$.
 
Let us now define for real numbers $x$ the function
\begin{equation}
S(x) = x^{\frac{1}{6}}\, \frac{\sum_{\mathfrak{p},\ N(\mathfrak{p})\le x} 
\dim H^1(\Gamma^0(\mathfrak{p}),\C)}
{|\{\mathfrak{p} \ \vrule\ N(\mathfrak{p})\le x\, \}|},
\end{equation}
where the sum is extended over all degree one prime ideals of
$\O_{-1}$. The function $S(x)$ can be tabulated in the range $x \le 60000$ as follows:
\begin{center}
\begin{tabular}{|c|c|c|c|c|c|c|c|c|c|c|}
\hline
$x/1000$  & $6$ & $12$ & $18$ & $24$ & $30$ & $36$ & $42$ & $48$ & $54$ & $60$\\
\hline
$S(x)$ & $3.21$ & $3.39$ & $3.62$ & $3.99$ & $4.15$ & $4.18$ & $4.24$ & $4.31$ & $4.37$ & $4.52$\\
\hline
\end{tabular}
\end{center}
See Question \ref{q13} of the introduction for some comments on this table.

\subsubsection{The case of non-trivial coefficients}\label{Gano}

We now report on some computational results on the cohomology
spaces $H^1(\Gamma^0(\mathfrak{p}),E_n)$ where $\mathfrak{p}$ is a prime of $\O_{-1}$ of degree one and $n \ge 1$.

Let $\pi$ be a generator of $\mathfrak{p}$ and let $\delta_\pi\in\GL(2,\Q(\sqrt{-1}))$ 
be defined as in (\ref{hedef}). The two injective homomorphisms
$$\iota_1 : \Gamma^0(\mathfrak{p})\to \SL(2,\O_{-1}), 
\qquad
\iota_2 : \Gamma^0(\mathfrak{p})\to \SL(2,\O_{-1}),
$$
where $\iota_1$ is just the injection and $\iota_2$ is induced by 
conjugation with the element $\delta_\pi$, give rise to an injection 
$$\iota: H^1(\SL(2,\O_{-1}),E_n)\oplus H^1(\SL(2,\O_{-1}),E_n)\hookrightarrow
H^1(\Gamma^0(\mathfrak{p}),E_n).$$
The image of $\iota$ is traditionally called the space of old classes. 
It is invariant under the Hecke operators and 
has an invariant complement $H_{\rm new}^1(\Gamma^0(\mathfrak{p}),E_n)$.
We have found: 
\begin{itemize}
\item $H_{\rm new}^1(\Gamma^0(\mathfrak{p}),E_1)=0$ for all prime ideals $\mathfrak{p}$ with
$N(\mathfrak{p})\le 1000$ except for the case $N(\mathfrak{p})=41$, where
$H_{\rm new}^1(\Gamma^0(\mathfrak{p}),E_1)$ has dimension $2$.
\item $H_{\rm new}^1(\Gamma^0(\mathfrak{p}),E_2)=0$ for all prime ideals $\mathfrak{p}$ with
$N(\mathfrak{p})\le 600$.
\item $H_{\rm new}^1(\Gamma^0(\mathfrak{p}),E_n)=0$ for all prime ideals $\mathfrak{p}$ with
$N(\mathfrak{p})\le 90$ and $3\le n\le 10$.
\end{itemize}
Table \ref{tahep2} contains some examples of the characteristic polynomials
of the Hecke operators on $H_{\rm new}^1(\Gamma^0(\mathfrak{p}),E_1)$ for $N(\mathfrak{p})=41$.
The Hecke operators on
$H_{\rm new}^1(\Gamma^0(\mathfrak{p}),E_1)$ satisfy
$T_{i\pi}=-T_\pi$ for all prime elements $\pi$ of degree one. There is no apparent 
connection between $T_\pi$ and $T_{\bar \pi}$.
We thank Haluk Sengun for help with this computation. 
\begin{table}
\begin{center}
\begin{tabular}{|c|c|c||c|c|c|}
\hline
$p$  & $\pi$ & $A_\pi$ & $p$ & $\pi$ & $A_\pi$ 
\\
\hline
$2$ & $1+i$ & $x^2 - x - 10$ & $-$ & $-$ & $-$\\
\hline
$5$ & $2+i$ & $(x+4)^2$ & $5$ & $2-i$ & $x^2+6x-32$\\
\hline
$13$ & $3+2i$ & $x^2+2x-40$ & $13$ & $3-2i$ & $(x+10)^2$\\
\hline
$17$ & $1+4i$ & $x^2-22x+80$ & $17$ & $1-4i$ & $x^2 + 24x -20$\\ 
\hline
$29$ & $5+2i$ & $x^2+48x-80$  & $29$ & $5-2i$ & $x^2-164$\\
\hline
$37$ & $1+6i$ & $x^2+4x-160$  & $37$ & $1-6i$ &  \\
\hline
$41$ & $5+4i$ & $x^2 + 48x + 412$ & $41$ & $5-4i$ & \\
\hline
$61$ & $6+5i$ & $x^2 -108+292$ & $61$ & $6-5i$ & \\
\hline
$73$ & $8+3i$ & $x^2 - 106x + 2440$ & $73$ & $8-3i$ & \\
\hline
\end{tabular}
\caption{\it Characteristic Polynomials of Hecke operators on 
$H_{\rm new}^1(\Gamma^0(\mathfrak{p}_{41}),E_1)$}\label{tahep2}
\end{center}
\end{table}

\section{Cohomology of non-arithmetic groups}\label{nonar}

This section contains computational results on 
the cohomology of various geometrically constructed and mostly
non-arithmetic groups. The results are 
discussed in more detail in the introduction. See Section \ref{Haeins} for remarks on the
method of computation and especially for the notation $\dim_{\le x}$ used below.

\subsection{Klimenko's examples}\label{lena}

The discrete subgroups $\G \subseteq \PSL(2,\C)$ described here arose in 
an important attempt to classify
simultaneous conjugacy classes of pairs of matrices 
generating discrete subgroups of $\SL(2,\C)$ (see
\cite{Klisa,Kliko1,Kliko2,Kliko3,Kliko4}). 
We follow the notation of
\cite{Kliko3,Kliko4}, see also \cite{Grukli}.

\subsubsection{Groups of finite covolume}

Let $k \ge 8$ be an even integer. We set
\begin{equation}\label{tee}
t=t_k=(\exp({\pi i}/{k})+\exp({-\pi i}/{k}))^2
=\exp({2\pi i}/{k})+\exp({-2\pi i}/{k})+2
\end{equation}
and define the matrices
\begin{equation}
f=f_k = \begin{pmatrix} \exp({\pi i}/{k}) & 0\\
0 & \exp(-{\pi i}/{k})
\end{pmatrix},
\end{equation}
\begin{equation}
g=g_k=\begin{pmatrix}
\frac{1}{2}\left( \sqrt{\frac{t}{(t-3)(4-t)}}+\sqrt{\frac{3}{t-3}}\right) & 1\\
\frac{t-3}{4-t} & 
\frac{1}{2}\left( \sqrt{\frac{t}{(t-3)(4-t)}}-\sqrt{\frac{3}{t-3}}\right)
\end{pmatrix}.
\end{equation}
Let 
${\rm GTet}_1[k,3,3]$ be the image in $\PSL (2,\C)$ of the
group generated by the matrices $f$ and $g$.
The following properties are known. 
\begin{itemize}
\item ${\rm GTet}_1[k,3,3]$ is a discrete  
subgroup of $\PSL(2,\C)$ of finite covolume with one cusp 
(\cite{Kliko3,Kliko4}).
\item ${\rm GTet}_1[k,3,3]$ is commensurable with a reflection group
(\cite{Kliko3,Kliko4}).
\item ${\rm GTet}_1[k,3,3]$ is non-arithmetic for all $k$ (\cite{Grukli}).
\item ${\rm GTet}_1[k,3,3]=\langle \, f,\, g\ \vrule\ f^k,\, 
(gf^{k/2}zf^{k/2}g^{-1}z)^3,\, z^2,\, fzgf^{-1}g^{-1}z\,\rangle$,
where $z=fgfg^{-1}f$ (\cite{Kliko3}).
\end{itemize}
Using this explicit presentation, we find:
\begin{equation}
\dim_{\le 1000} H^1({\rm GTet}_1[k,3,3], E_n)=
\begin{cases}
2(n-1)/3+1,\qquad \qquad \qquad \ {\rm if}\ n\equiv 1 \ (6),\\
2(n-2)/3+2,\qquad \qquad \qquad \ {\rm if}\ n\equiv 2 \ (6),\\
2(n-3)/3+2,\qquad \qquad \qquad \ {\rm if}\ n\equiv 3 \ (6),\\
2(n-4)/3+3,\qquad \qquad \qquad \ {\rm if}\ n\equiv 4 \ (6),\\
2(n-5)/3+4,\qquad \qquad \qquad \ {\rm if}\ n\equiv 5 \ (6),\\ 
2(n-6)/3+4,\qquad \qquad \qquad \ {\rm if}\ n\equiv 6 \ (6),\\ 
\end{cases}
\end{equation}
in the range $8\le k\le 100$ and $1\le n \le 50$. Note that these groups $\Gamma$ are invariant under the complex conjugation automorphism of $\PSL (2,\C)$.
This opens up the possibility to compute the trace of this involution 
on $H^1_{\rm cusp} (\Gamma, E_n)$ and to obtain a lower bound for the dimension of this space, following the work of Rohlfs and Kr\"amer
\cite{Kraemer,Rohlfs} on the Bianchi groups. We hope to come back to this question in the future.

\subsubsection{Cocompact groups}

As in the case considered before, we take from 
 \cite{Kliko3} (see also \cite{Grukli}) a
series of explicit pairs of matrices generating a discrete subgroup 
$\G\subseteq \PSL(2,\C)$. In this case the groups $\Gamma$ act on 
three-dimensional hyperbolic space with a compact quotient. 

Let $k\ge 8$ be an even integer. We define 
$t=t_k$ as in (\ref{tee}) and set
\begin{equation}
f=f_k=\begin{pmatrix} \exp({\pi i}/{k}) & 0\\
0 & \exp(-{\pi i}/{k})
\end{pmatrix},
\end{equation}
\begin{equation}
g=g_k=\begin{pmatrix}
\frac{1}{2}\left( \sqrt{\frac{2(t-2)}{(t-3)(4-t)}}
+\sqrt{\frac{2}{t-3}}\right) & 1\\
\frac{t-3}{4-t} & 
\frac{1}{2}\left( \sqrt{\frac{2(t-2)}{(t-3)(4-t)}}-\sqrt{\frac{2}{t-3}}\right)
\end{pmatrix}.
\end{equation}
Define ${\rm GTet}_1[k,3,2]\subseteq \PSL(2,\C)$ to be the image in 
$\PSL (2,\C)$ of the group generated by
$f$ and $g$. The following facts are known.
\begin{itemize}
\item ${\rm GTet}_1[k,3,2]$ is a discrete and cocompact subgroup 
of $\PSL(2,\C)$ (\cite{Kliko3,Kliko4}). 
\item ${\rm GTet}_1[k,3,2]$ is commensurable with a reflection group
(\cite{Kliko3,Kliko4}).
\item ${\rm GTet}_1[k,3,2]$ is non-arithmetic for all $k\ge 14$ 
(\cite{Grukli}).
\item ${\rm GTet}_1[k,3,2]=\langle \, f,\, g\ \vrule\ f^k,\, 
(gf^{k/2}zf^{k/2}g^{-1}z)^2,\, z^2,\, fzgf^{-1}g^{-1}z\,\rangle$,
where $z=fgfg^{-1}f$ (\cite{Kliko3}).
\end{itemize}
Again these results are sufficient to compute cohomology spaces. We find
\begin{equation}\label{coho5}
\dim_{\le 1000} H^1({\rm GTet}_1[k,3,2], E_n)=
\begin{cases}
n/4, \qquad\qquad\qquad\qquad  \quad {\rm if}\ n\equiv 0 \ (4), \\
(n+1)/2,\qquad \qquad \qquad \ {\rm if}\ n\equiv 1 \ (2),\\
(n+2)/4, \qquad \qquad\qquad \ {\rm if}\ n\equiv 2 \ (4)
\end{cases}
\end{equation}
in the range $14\le k\le 100$ and $1\le n \le 30$. The groups 
${\rm GTet}_1[8,3,2]$, ${\rm GTet}_1[10,3,2]$ and ${\rm GTet}_1[12,3,2]$ 
are arithmetic. Compared to (\ref{coho5}), the dimensions of their cohomology groups show a similar but slightly more complicated behavior. In particular, the dimensions 
of the cohomology spaces $H^1( {\rm GTet}_1[10,3,2],E_n)$ are given by linear functions on the residue classes modulo $20$ within the range of our computations.
Again all these groups are invariant under the complex conjugation automorphism of $\PSL (2,\C)$.

\subsection{Helling's examples}\label{hell}

Here we report on a series of two-generator discrete subgroups of
$\SL(2,\C)$ described in \cite{hell}. We shall keep the terminology of 
\cite{hell}. 
The phenomena seen here are new.

For a non-negative integer $k$, let $T_k$ and $U_k$ be the standard Tchebyshev 
polynomials \cite{Mag}. They can be defined by the relation
\[
\left( x + \sqrt{x^2-1} \right)^k = T_k (x) + U_{k-1} (x) \sqrt{x^2-1},
\]
for example.
For a non-negative integer $m$ we define
polynomials 
\begin{equation}
\tilde{p}_m(x)=
        \begin{cases} 2T_k\left(\frac{x}{2}\right),\qquad \qquad\qquad\ \ 
{\rm if}\ m=2k,\\
  U_k\left(\frac{x}{2}\right)-U_{k-1}\left(\frac{x}{2}\right),\qquad 
{\rm if}\ m=2k+1,              
        \end{cases} 
\end{equation}
and 
\begin{equation}
{f}_m(x)=
        \begin{cases} \tilde{p}_{m+2}(x)^2-x^2+4, \qquad
{\rm if}\ m \ {\rm is\ even},\\
\tilde{p}_{m+2}(x)^2-x+2, \qquad \ \ 
           {\rm if}\ m \ {\rm is\ odd.}
        \end{cases} 
\end{equation}

The following table contains the first ten polynomials $f_m(x)$.

\bigskip
\begin{center}
\begin{tabular}{|c|c|}
\hline
$m$  & $f_m(x)$ \\
\hline
$1$ & $x^2 - 3x + 3$\\
$2$ & $x^4 - 5x^2 + 8$\\
$3$ & $x^4 - 2x^3 - x^2 + x + 3$\\
$4$ & $x^6 - 6x^4 + 8x^2 + 4$\\
$5$ & $x^6 - 2x^5 - 3x^4 + 6x^3 + 2x^2 - 5x + 3$\\
$6$ & $x^8 - 8x^6 + 20 x^4 - 17x^2 + 8$\\
$7$ & $x^8 - 2x^7 - 5x^6 + 10x^5 + 7x^4 - 14x^3 - 2x^2 + 3x + 3$\\
$8$ & $x^{10} - 10x^8 + 35x^6 - 50x^4 + 24x^2 + 4$\\
$9$ & $x^{10} - 2x^9- 7x^8+ 14x^7+ 16x^6 - 32x^5-13x^4+26x^3 + 3x^2-7x+3$\\
$10$ & $x^{12} - 12x^{10} + 54x^8 - 112x^6 + 105x^4 - 37x^2 + 8$\\
\hline
\end{tabular}
\end{center}
Helling shows in \cite{hell} that the polynomials $f_m$ 
have 
only non-real zeroes. For a zero $z$ of $f_m$ define the 
matrices
\begin{equation}
A_m=\begin{pmatrix} 0 & 1\\
-1 & z
\end{pmatrix},\quad 
B_m=\begin{pmatrix} 1 & 0\\
\frac{\tilde{p}_m(z)}{\tilde{p}_{m+2}(z)} & 1
\end{pmatrix},\quad 
C_m=\begin{pmatrix} 1 & \frac{\tilde{p}_m(z)}{\tilde{p}_{m+2}(z)}\\ 
0 & 1
\end{pmatrix}.
\end{equation}

An easy computation using properties of the Tchebyshev polynomials confirms 
that these matrices satisfy the relations
\begin{equation}\label{rela5}
A_mC_mA_m^{-1}=B_m^{-1},\qquad C_mB_mC_m^{-1}B_m^{-1}=A_m^m.
\end{equation}
Define $\Theta_m$ to be the group generated by the above matrices:
\begin{equation}
\Theta_m=\langle A_m,\, B_m, C_m\rangle=\langle A_m,\,  C_m\rangle
\subseteq \SL(2,\C).
\end{equation}
Helling shows that for every $m\in \mathbb{N}$ there is a zero 
$z\in \C$ of $f_m$ such that the matrix group $\Theta_m$ satisfies:
\begin{itemize}
\item $\Theta_m$ is a discrete and torsion-free subgroup of $\SL(2,\C)$.
\item $\Theta_m$ has finite covolume and exactly one cusp.
\item $\Theta_m$ is defined by the relations (\ref{rela5}).
\item The groups $\Theta_1$, $\Theta_2$ are arithmetic groups, but all the
other $\Theta_m$ ($m\ge 3$) are non-arithmetic. 
\end{itemize}
The zero $z$ in question is specified (up to complex conjugation) by the condition
\[
|z-2| < 4 \sin^2 \frac{\pi}{2m}
\]
for $m \ge 3$ odd and
\[
{\rm Re} \, (z) > 0, \quad  |z^2 - 4| < 4 \sin^2 \frac{\pi}{m}
\]
for $m \ge 4$ even. For $m = 1$ or $2$, $z$ may be taken to be any zero of $f_m$. 
Concerning the cohomology of the groups $\Theta_m$, we 
can report the following computations:

\bigskip
\centerline{$\bf m=1$:}
The group $\Theta_1$ is (up to conjugacy) the famous figure eight knot group.
It is conjugate to a congruence subgroup of $\SL(2,{\cal O}_{-3})$.
For $k\le 120$ we have

\begin{equation}
\dim H^1(\Theta_1, E_n)=
\begin{cases}
n/3,\qquad \qquad \qquad \ {\rm if}\ n\equiv 0 \ {\rm mod}\ 3,\\
(n+2)/3, \qquad \qquad {\rm if}\ n\equiv 1 \ {\rm mod}\ 3, \\
(n+1)/3+1, \qquad \ {\rm if}\ n\equiv 2 \ {\rm mod}\ 3. 
\end{cases}
\end{equation}

\bigskip
\centerline{$\bf m=2$:}

The group $\Theta_2$ is isomorphic to the fundamental group of the 
lens space with fundamental group of order $2$
with a knot removed.  
It is conjugate to a group commensurable with $\SL(2,{\cal O}_{-7})$.
For $k\le 120$ we have
\begin{equation}
\dim H^1(\Theta_2, E_n)=
\begin{cases}
n/3,\qquad \qquad \qquad \ {\rm if}\ n\equiv 0 \ {\rm mod}\ 3,\\
(n+2)/3, \qquad \qquad {\rm if}\ n\equiv 1 \ {\rm mod}\ 3, \\
(n+1)/3,\qquad \qquad {\rm if}\ n\equiv 2 \ {\rm mod}\ 3, 
\end{cases}
\end{equation}
except in the case $n=12$, where we have
$$\dim H^1(\Theta_2, E_{12})=6.$$

\bigskip
\centerline{$\bf m\ge 3$:}

Here we have 
\begin{equation}
\dim H^1(\Theta_m, E_{n})=1
\end{equation}
for all $3\le m\le 150$ and $1\le k\le 30$. This means that in this range we have $H_{\rm cusp}^1(\Theta_m, E_{n})=0$.

\subsection{A cocompact tetrahedral group}\label{tetra}

Let ${\bf CT}(26)$ be the tetrahedral hyperbolic reflection group 
constructed (for example) in Section 10 of \cite{EGM3} and let 
$\Gamma_{26}\subseteq \PSL(2,\C)$ be its unique subgroup of index $2$. 
The quotient
$\PSL(2,\C)/ \Gamma_{26}$ is compact. 
The group  $\Gamma_{26}$ is non-arithmetic and has the 
presentation
\begin{equation} 
\Gamma_{26}=\langle\,  a,b,c \ \vrule \ a^3,\, b^2,\, c^5,\,
(ac^{-1})^2,\, (bc^{-1})^3,\, (ab)^4\, \rangle.
\end{equation}
Using the data from \cite{EGM3} we infer that  $\Gamma_{26}$ can be generated 
(up to conjugacy) by the matrices 
$$
a=\begin{pmatrix}  \frac{2t^3 + t^2 + t + 2}{5}
 & 1\\
\frac{-t^3 +t^2 - 2}{5}      & \frac{-t^3 -t^2 -t+3}{5}
\end{pmatrix}, \quad
b=\begin{pmatrix}  \frac{-3t^3 +t^2 - 4t + 2}{5}
 & b_2\\
c_2      & \frac{3t^3 - t^2 + 4t - 2}{5}
\end{pmatrix}, \quad
c=\begin{pmatrix} t^{-1} & 0\\
0      & t
\end{pmatrix},
$$
where $t\in \C$ is a primitive $10$-th root of unity and 
$c_2$ is one of the two complex roots of the polynomial
$$
x^4 + \frac{-6t^3 + 6t^2 + 8}{5} x^3 + \frac{-t^3 + t^2 - 3}{5}x^2 
+ \frac{-4t^3 + 4t^2 + 2}{25}x + \frac{3t^3 - 3t^2 + 2}{25}.
$$
The entry $b_2$ is determined by
$$
b_2=(-20t^3 + 20t^2 + 35)c_2^3 
+ (-50t^3 + 50t^2 + 80)c_2^2 + (9t^3 - 9t^2 -17)c_2 - 4t^3 + 4t^2 + 6.
$$
We have found that
\begin{equation}
H^1(\Gamma_{26},E_n)=0
\end{equation}
for $0\le n\le 90$.
The group $\Gamma_{26}$ has $222$ conjugacy classes of 
subgroups of index less than
or equal to $24$. We have also determined the 
dimensions of some cohomology spaces of these subgroups.
Of the $222$ subgroups $191$ satisfied
$H^1(\Gamma,E_n)=0$ 
in the range $0\le n\le 10$. Thirty subgroups had  
$\dim_{\le 1000} \, H^1(\Gamma,E_n)=1$ in the 
range $0\le n\le 10$. One of the $222$ had   
$\dim_{\le 1000} \, H^1(\Gamma,E_n)=2$, again in this range.


\end{document}